\documentclass[12pt]{amsart}
\usepackage{amsmath,amssymb,amsfonts}
\begin{document}

\theoremstyle{plain}
\newtheorem{thm}{Theorem}[subsection]
\newtheorem{lem}[thm]{Lemma}
\newtheorem{cor}[thm]{Corollary}
\newtheorem{prop}[thm]{Proposition}
\newtheorem{remark}[thm]{Remark}
\newtheorem{defn}[thm]{Definition}
\newtheorem{ex}[thm]{Example}

\numberwithin{equation}{subsection}
\newcommand{\mc}{\mathcal}
\newcommand{\mb}{\mathbb}
\newcommand{\surj}{\twoheadrightarrow}
\newcommand{\inj}{\hookrightarrow}
\newcommand{\red}{{\rm red}}
\newcommand{\codim}{{\rm codim}}
\newcommand{\rank}{{\rm rank}}
\newcommand{\Pic}{{\rm Pic}}
\newcommand{\Div}{{\rm Div}}
\newcommand{\Hom}{{\rm Hom}}
\newcommand{\im}{{\rm im}}
\newcommand{\Spec}{{\rm Spec \,}}
\newcommand{\Sing}{{\rm Sing}}
\newcommand{\Char}{{\rm char}}
\newcommand{\Tr}{{\rm Tr}}
\newcommand{\Gal}{{\rm Gal}}
\newcommand{\Min}{{\rm Min \ }}
\newcommand{\Max}{{\rm Max \ }}
\newcommand{\sA}{{\mathcal A}}
\newcommand{\sB}{{\mathcal B}}
\newcommand{\sC}{{\mathcal C}}
\newcommand{\sD}{{\mathcal D}}
\newcommand{\sE}{{\mathcal E}}
\newcommand{\sF}{{\mathcal F}}
\newcommand{\sG}{{\mathcal G}}
\newcommand{\sH}{{\mathcal H}}
\newcommand{\sI}{{\mathcal I}}
\newcommand{\sJ}{{\mathcal J}}
\newcommand{\sK}{{\mathcal K}}
\newcommand{\sL}{{\mathcal L}}
\newcommand{\sM}{{\mathcal M}}
\newcommand{\sN}{{\mathcal N}}
\newcommand{\sO}{{\mathcal O}}
\newcommand{\sP}{{\mathcal P}}
\newcommand{\sQ}{{\mathcal Q}}
\newcommand{\sR}{{\mathcal R}}
\newcommand{\sS}{{\mathcal S}}
\newcommand{\sT}{{\mathcal T}}
\newcommand{\sU}{{\mathcal U}}
\newcommand{\sV}{{\mathcal V}}
\newcommand{\sW}{{\mathcal W}}
\newcommand{\sX}{{\mathcal X}}
\newcommand{\sY}{{\mathcal Y}}
\newcommand{\sZ}{{\mathcal Z}}
\newcommand{\A}{{\Bbb A}}
\newcommand{\B}{{\Bbb B}}
\newcommand{\C}{{\Bbb C}}
\newcommand{\D}{{\Bbb D}}
\newcommand{\E}{{\Bbb E}}
\newcommand{\F}{{\Bbb F}}
\newcommand{\G}{{\Bbb G}}
\renewcommand{\H}{{\Bbb H}}
\newcommand{\I}{{\Bbb I}}
\newcommand{\J}{{\Bbb J}}
\newcommand{\M}{{\Bbb M}}
\newcommand{\N}{{\Bbb N}}
\renewcommand{\P}{{\Bbb P}}
\newcommand{\Q}{{\Bbb Q}}
\newcommand{\R}{{\Bbb R}}
\newcommand{\T}{{\Bbb T}}
\newcommand{\U}{{\Bbb U}}
\newcommand{\V}{{\Bbb V}}
\newcommand{\W}{{\Bbb W}}
\newcommand{\X}{{\Bbb X}}
\newcommand{\Y}{{\Bbb Y}}
\newcommand{\Z}{{\Bbb Z}}

\catcode`\@=11
%
%
\def\opn#1#2{\def#1{\mathop{\kern0pt\fam0#2}\nolimits}} 
\def\bold#1{{\bf #1}}%
\def\underrightarrow{\mathpalette\underrightarrow@}
\def\underrightarrow@#1#2{\vtop{\ialign{$##$\cr
 \hfil#1#2\hfil\cr\noalign{\nointerlineskip}%
 #1{-}\mkern-6mu\cleaders\hbox{$#1\mkern-2mu{-}\mkern-2mu$}\hfill
 \mkern-6mu{\to}\cr}}}
\let\underarrow\underrightarrow
\def\underleftarrow{\mathpalette\underleftarrow@}
\def\underleftarrow@#1#2{\vtop{\ialign{$##$\cr
 \hfil#1#2\hfil\cr\noalign{\nointerlineskip}#1{\leftarrow}\mkern-6mu
 \cleaders\hbox{$#1\mkern-2mu{-}\mkern-2mu$}\hfill
 \mkern-6mu{-}\cr}}}
%
%

%
\def\:{\colon}
\let\oldtilde=\tilde
\def\tilde#1{\mathchoice{\widetilde{#1}}{\widetilde{#1}}%
{\indextil{#1}}{\oldtilde{#1}}}
\def\indextil#1{\lower2pt\hbox{$\textstyle{\oldtilde{\raise2pt%
\hbox{$\scriptstyle{#1}$}}}$}}
\def\pnt{{\raise1.1pt\hbox{$\textstyle.$}}}
%
%
\let\amp@rs@nd@\relax
\newdimen\ex@
\ex@.2326ex
\newdimen\bigaw@
\newdimen\minaw@
\minaw@16.08739\ex@
\newdimen\minCDaw@
\minCDaw@2.5pc
\newif\ifCD@
\def\minCDarrowwidth#1{\minCDaw@#1}
\newenvironment{CD}{\@CD}{\@endCD}
\def\@CD{\def\A##1A##2A{\llap{$\vcenter{\hbox
 {$\scriptstyle##1$}}$}\Big\uparrow\rlap{%
$\vcenter{\hbox{$\scriptstyle##2$}}$}&&}%
\def\V##1V##2V{\llap{$\vcenter{\hbox
 {$\scriptstyle##1$}}$}\Big\downarrow\rlap{$\vcenter{%
\hbox{$\scriptstyle##2$}}$}&&}%
\def\={&\hskip.5em\mathrel
 {\vbox{\hrule width\minCDaw@\vskip3\ex@\hrule width
 \minCDaw@}}\hskip.5em&}%
\def\verteq{\Big\Vert&&}%
\def\noarr{&&}%
\def\vspace##1{\noalign{\vskip##1\relax}}\relax\let\amp@rs@nd@&\iffalse}\fi \CD@true\vcenter\bgroup\relax\let\\=\cr\iffalse}\fi\tabskip\z@skip\baselineskip20\ex@
 \lineskip3\ex@\lineskiplimit3\ex@\halign\bgroup
 &\hfill$\m@th##$\hfill\cr}
\def\@endCD{\cr\egroup\egroup}
%
\def\>#1>#2>{\amp@rs@nd@\setbox\z@\hbox{$\scriptstyle
 \;{#1}\;\;$}\setbox\@ne\hbox{$\scriptstyle\;{#2}\;\;$}\setbox\tw@
 \hbox{$#2$}\ifCD@
 \global\bigaw@\minCDaw@\else\global\bigaw@\minaw@\fi
 \ifdim\wd\z@>\bigaw@\global\bigaw@\wd\z@\fi
 \ifdim\wd\@ne>\bigaw@\global\bigaw@\wd\@ne\fi
 \ifCD@\hskip.5em\fi
 \ifdim\wd\tw@>\z@
 \mathrel{\mathop{\hbox to\bigaw@{\rightarrowfill}}\limits^{#1}_{#2}}\else
 \mathrel{\mathop{\hbox to\bigaw@{\rightarrowfill}}\limits^{#1}}\fi
 \ifCD@\hskip.5em\fi\amp@rs@nd@}
\def\<#1<#2<{\amp@rs@nd@\setbox\z@\hbox{$\scriptstyle
 \;\;{#1}\;$}\setbox\@ne\hbox{$\scriptstyle\;\;{#2}\;$}\setbox\tw@
 \hbox{$#2$}\ifCD@
 \global\bigaw@\minCDaw@\else\global\bigaw@\minaw@\fi
 \ifdim\wd\z@>\bigaw@\global\bigaw@\wd\z@\fi
 \ifdim\wd\@ne>\bigaw@\global\bigaw@\wd\@ne\fi
 \ifCD@\hskip.5em\fi
 \ifdim\wd\tw@>\z@
 \mathrel{\mathop{\hbox to\bigaw@{\leftarrowfill}}\limits^{#1}_{#2}}\else
 \mathrel{\mathop{\hbox to\bigaw@{\leftarrowfill}}\limits^{#1}}\fi
 \ifCD@\hskip.5em\fi\amp@rs@nd@}
%
%
\newenvironment{CDS}{\@CDS}{\@endCDS}
\def\@CDS{\def\A##1A##2A{\llap{$\vcenter{\hbox
 {$\scriptstyle##1$}}$}\Big\uparrow\rlap{$\vcenter{\hbox{%
$\scriptstyle##2$}}$}&}%
\def\V##1V##2V{\llap{$\vcenter{\hbox
 {$\scriptstyle##1$}}$}\Big\downarrow\rlap{$\vcenter{\hbox{%
$\scriptstyle##2$}}$}&}%
\def\={&\hskip.5em\mathrel
 {\vbox{\hrule width\minCDaw@\vskip3\ex@\hrule width
 \minCDaw@}}\hskip.5em&}
\def\verteq{\Big\Vert&}
\def\novarr{&}
\def\noharr{&&}
\def\SE##1E##2E{\slantedarrow(0,18)(4,-3){##1}{##2}&}
\def\SW##1W##2W{\slantedarrow(24,18)(-4,-3){##1}{##2}&}
\def\NE##1E##2E{\slantedarrow(0,0)(4,3){##1}{##2}&}
\def\NW##1W##2W{\slantedarrow(24,0)(-4,3){##1}{##2}&}
\def\slantedarrow(##1)(##2)##3##4{\thinlines\unitlength1pt\lower 6.5pt\hbox{%
\begin{picture}(24,18)%
\put(##1){\vector(##2){24}}%
\put(0,8){$\scriptstyle##3$}%
\put(20,8){$\scriptstyle##4$}%
\end{picture}}}
\def\vspace##1{\noalign{\vskip##1\relax}}\relax\let\amp@rs@nd@&\iffalse}\fi
 \CD@true\vcenter\bgroup\relax\let\\=\cr\iffalse}\fi\tabskip\z@skip\baselineskip20\ex@
 \lineskip3\ex@\lineskiplimit3\ex@\halign\bgroup
 &\hfill$\m@th##$\hfill\cr}
\def\@endCDS{\cr\egroup\egroup}
%
\newdimen\TriCDarrw@
\newif\ifTriV@
\newenvironment{TriCDV}{\@TriCDV}{\@endTriCD}
\newenvironment{TriCDA}{\@TriCDA}{\@endTriCD}
\def\@TriCDV{\TriV@true\def\TriCDpos@{6}\@TriCD}
\def\@TriCDA{\TriV@false\def\TriCDpos@{10}\@TriCD}
\def\@TriCD#1#2#3#4#5#6{%
\setbox0\hbox{$\ifTriV@#6\else#1\fi$}
\TriCDarrw@=\wd0 \advance\TriCDarrw@ 24pt
\advance\TriCDarrw@ -1em
\def\SE##1E##2E{\slantedarrow(0,18)(2,-3){##1}{##2}&}
\def\SW##1W##2W{\slantedarrow(12,18)(-2,-3){##1}{##2}&}
\def\NE##1E##2E{\slantedarrow(0,0)(2,3){##1}{##2}&}
\def\NW##1W##2W{\slantedarrow(12,0)(-2,3){##1}{##2}&}
\def\slantedarrow(##1)(##2)##3##4{\thinlines\unitlength1pt
\lower 6.5pt\hbox{\begin{picture}(12,18)%
\put(##1){\vector(##2){12}}%
\put(-4,\TriCDpos@){$\scriptstyle##3$}%
\put(12,\TriCDpos@){$\scriptstyle##4$}%
\end{picture}}}
\def\={\mathrel {\vbox{\hrule
   width\TriCDarrw@\vskip3\ex@\hrule width
   \TriCDarrw@}}}
\def\>##1>>{\setbox\z@\hbox{$\scriptstyle
 \;{##1}\;\;$}\global\bigaw@\TriCDarrw@
 \ifdim\wd\z@>\bigaw@\global\bigaw@\wd\z@\fi
 \hskip.5em
 \mathrel{\mathop{\hbox to \TriCDarrw@
{\rightarrowfill}}\limits^{##1}}
 \hskip.5em}
\def\<##1<<{\setbox\z@\hbox{$\scriptstyle
 \;{##1}\;\;$}\global\bigaw@\TriCDarrw@
 \ifdim\wd\z@>\bigaw@\global\bigaw@\wd\z@\fi
 \mathrel{\mathop{\hbox to\bigaw@{\leftarrowfill}}\limits^{##1}}
 }
 \CD@true\vcenter\bgroup\relax\let\\=\cr\iffalse}\fi
 \tabskip\z@skip\baselineskip20\ex@
 \lineskip3\ex@\lineskiplimit3\ex@
 \ifTriV@
 \halign\bgroup
 &\hfill$\m@th##$\hfill\cr
#1&\multispan3\hfill$#2$\hfill&#3\\
&#4&#5\\
&&#6\cr\egroup%
\else
 \halign\bgroup
 &\hfill$\m@th##$\hfill\cr
&&#1\\%
&#2&#3\\
#4&\multispan3\hfill$#5$\hfill&#6\cr\egroup
\fi}
\def\@endTriCD{\egroup}


\title{Algebraic Chern-Simons Theory}
\author{Spencer Bloch}
\address{Department of Mathematics, University of Chicago, Chicago, IL 60637}
\email{bloch@math.uchicago.edu}
\thanks{The first author is supported in part by the NSF, the second one is
supported in part by the DFG}
\author{H\'el\`ene Esnault}
\address{Universit\"at Essen,FB 6, Mathematik, 45117 Essen, Germany}
\email{esnault@uni-essen.de}
\maketitle
\setcounter{section}{-1}

\section{Introduction}\label{sec:intro}
\subsection{}Secondary (Chern-Simons) characteristic classes associated to
bundles with connection play an important role in differential geometry. We
propose to investigate a related construction for algebraic bundles. Non-flat
algebraic connections for bundles on complex projective manifolds are virtually
non-existent (we know of none), and a deep theorem of Reznikov \cite{Re}
implies
that Chern-Simons classes are torsion for flat bundles on such spaces. On the
other hand, it is possible (in several different ways, cf.
[\ref{subsec:aff_concl}] below)
given a vector bundle $E$ on $X$ to construct an affine fibration $f : Y \to X$
(i.e.
locally over $X$,
$Y\cong X\times \A^n$) such that $f^*E$ admits an algebraic connection. One can
arrange moreover that $Y$ itself be an affine variety. Since
pullback $f^*$ induces an isomorphism from the {\it chow motive} of $X$ to that
of $Y$, one can in some sense say that every algebraic variety is
equivalent to an affine variety, and every vector bundle is equivalent to a
vector bundle with an algebraic connection. Thus, an algebraic Chern-Simons
theory has some interest. Speaking loosely, the content of such a theory is
that a
closed differential form $\tau$ representing a characteristic class like the
chern class of a vector
bundle on a variety $X$ will be Zariski-locally exact, $\tau|U_i = d\eta_i$.
The choice of a
connection on the bundle enables one to choose the primitives $\eta_i$
canonically upto an
exact form. In particular, $(\eta_i-\eta_j)|U_i\cap U_j$ is exact. When $X$ is
affine, a different
choice of connection will change the $\eta_i$ by a global form $\eta$.

\subsection{}Unless otherwise noted, all our spaces $X$ will be smooth,
quasi-projective varieties over a field $k$ of characteristic $0$. Given a
bundle of rank $N$ with connection $(E,\nabla)$ on $X$ and an invariant
polynomial
$P$ of degree $n$ on the Lie algebra of ${\rm GL}_N$ (cf. \cite{CS}), we
construct
classes

\begin{equation}
w_n(E,\nabla,P)\in \Gamma(X,\Omega_X^{2n-1}/d\Omega^{2n-2}_X);\ n\ge 2.
\end{equation}
Here $\Omega^i_X$ is the Zariski sheaf of K\" ahler $i$-forms on
$X$, and
$d :\Omega^i_X \to \Omega^{i+1}_X$ is exterior differentiation. Zariski
locally,
these classes are given explicitly in terms of universal polynomials in the
connection and its curvature. They satisfy the basic compatibility:
\begin{quote}
$dw_n(E,\nabla,P)$ is a closed $2n$-form representing the characteristic class
in de Rham cohomology associated to $P$ by Chern-Weil theory. Note that $dw_n$
is not
necessarily exact, because $w_n$ is not a globally defined form.
\end{quote}

The simplest example is to take $E$ trivial of rank $2$ and to assume the
connection on the
determinant bundle is trivial. The connection is then given by a matrix of
$1$-forms
$A=\bigl(\begin{smallmatrix} \alpha & \beta\\\gamma & -\alpha
\end{smallmatrix}\bigr)$.
Taking
$P(M)=\text{Tr}(M^2)$ one finds
\begin{gather}\label{eqn:rk2ex}
w_2(E,\nabla,P) = 2\alpha\wedge d\alpha -4 \alpha \beta \gamma + \beta d\gamma
+ \gamma d\beta \notag \\
\mbox{ or, if } A \mbox{ is integrable, } \\
w_2(E, \nabla, P) = -2 \alpha \wedge d\alpha = -2 \alpha \beta \gamma. \notag
\end{gather}

One particularly important invariant polynomial $P_n$ maps a diagonal matrix to
the $n$-th
elementary symmetric function in its entries. We write
\begin{equation}\label{eqn0.2.2} w_n(E,\nabla) := w_n(E,\nabla,P_n)
\end{equation}
For example, $P_2(M) := \frac{1}{2}({\rm Tr}M)^2 - {\rm
Tr}(M^2)$. In fact, when $\nabla$ is integrable, $w_n(E, \nabla,
P) = \lambda w_n(E, \nabla)$ for some coefficient $\lambda \in
\Q$ (see \ref{prop:unicity}).

When $k=\C$, $w_n(E,\nabla)$ is linked to the chern class in $A^n(X)$, where
$A^n(X)$ denotes the group of algebraic cycles modulo a certain adequate
equivalence relation,
homological equivalence on a divisor. For example, $A^2(X)$ is the group of
codimension $2$ cycles
modulo algebraic equivalence. When $n=2$ and $X$ is affine, there is an
isomorphism
\begin{equation}
\varphi :
\Gamma(X,\Omega_X^{3}/d\Omega^{2}_X)/\Gamma(X,\Omega_X^{3})\cong
A^2(X)\otimes_\Z \C.
\end{equation}

Writing $c_{2,{\rm cycle}}(E)$ for the second chern class of $E$ in $A^2(X)$,
we
have
\begin{equation} \varphi(w_2(E,\nabla)) = c_{2,{\rm cycle}}(E)\otimes 1
\end{equation}

\subsection{} Suppose now the connection $\nabla$ on $E$ is integrable, i.e.
$E$ is flat. Let $\sK_i^m$ denote the Zariski sheaf, image of
the Zariski Milnor $K$ sheaf in the constant sheaf $K^M_i(k(X))$.
One has a map ${\rm dlog} : \sK_i^m \to \Omega^i_{X,{\rm
clsd}}$. Functorial and additive classes
\begin{equation}\label{eqn0.3.1} c_i(E,\nabla) \in \H^i(X,\sK_i^m\to
\Omega^i \to \Omega^{i+1}\to
\ldots)
\end{equation}
were constructed in \cite{EII}. One has a natural map of complexes
\begin{equation}
\sigma : \{\sK_i^m\to \Omega^i \to \Omega^{i+1}\to\ldots\} \to
\Omega_X^{2i-1}/d\Omega^{2i-2}_X[-i].
\end{equation}

We prove in section \ref{sec:thetaw}
\begin{equation}\label{eqn0.3.3}
w_i(E,\nabla) = \sigma(c_i(E,\nabla))\in
\Gamma(X,\Omega^{2i-1}/d\Omega^{2i-2}).
\end{equation}
In the case of an integrable connection, the classes $w_n(E,\nabla)$ are
closed.
We are unable to answer the following
\subsubsection{basic question.} Are the classes $w_i(E,\nabla,P)$ all zero for
an
integrable connection $\nabla$?\label{subsubsec:question}

\subsection{}We continue to assume $\nabla$ integrable. We take $k=\C$, and $X$
smooth and
projective. We define the (generalized) Griffiths group $\text{Griff}^{n}(X)$
to be the group of
algebraic cycles of codimension $n$ homologous to zero, modulo those homologous
to zero on a
divisor. (For $n=2$, this is the usual Griffiths group of codimension $2$
algebraic cycles
homologous to zero modulo algebraic equivalence.) Our main result is
\begin{thm} We have $w_n(E,\nabla)=0$ if and only if $c_n(E)=0$ in
$\text{Griff}^n(X)\otimes\mb Q$.
\end{thm}
The proof of this theorem is given in section \ref{sec:griff}.

The idea is that one can associate to
any codimension $n$ cycle $Z$ homologous to zero
an extension of mixed Hodge structures of $\mb
Q(0)$ by $H^{2n-1}(X,\mb Q(n))$. One gets a quotient extension
$$0 \to H^{2n-1}(X,\mb Q(n))/N^1 \to E \to
\text{Griff}^n(X)\otimes\mb Q(0) \to 0
$$
where $N^1$ is the subspace of ``coniveau'' $1$,
the group on the right has the trivial Hodge
structure and where $$E\subset H^0(X,\mc
H^{2n-1}(\mb Q(n))) .$$

Using the classes \eqref{eqn0.3.1} and the comparison \eqref{eqn0.3.3} we show
$$w_n(E,\nabla)\in F^0E\cap E(\mathbb R).
$$
Furthermore,
$w_n(E,\nabla)\in E(\mb C)$ maps to the class of $c_n(E)$.
Since the kernel of this extension is pure of weight
$-1$ it follows easily that $w_n = 0 \Leftrightarrow c_n=0$.
In fact, Reznikov's theorem \cite{Re} implies
$$w_n(E,\nabla)\in E(\mb Q).
$$

\subsection{}
Through its link to the Griffiths group, it is clear that
the classes $w_n(E, \nabla)$, when $\nabla$ is integrable, are
rigid in a variation of the flat bundle $(E, \nabla)$ over
$X$. But in fact, a stronger rigidity ( see \ref{prop:rigid})
holds true: one can allow a 1 dimensional variation of
$X$ as well.

\subsection{} Examples (including Gau{\ss}-Manin systems of semi-stable
families of curves,
weight 1 Gau{\ss}-Manin systems, weight 2 Gau{\ss}-Manin systems of surfaces,
and
local systems with finite monodromy) for which the classes
$w_n(E, \nabla)$  vanish
are discussed in section \ref{sec:ex}.

It is possible (cf. section \ref{sec:ex}) to
define $w_n(E,\nabla,P)$ in characteristic $p$. In arithmetic situations, the
resulting classes are compatible with reduction mod $p$. When the bundle
$(E,\nabla)$ in characteristic $p$ comes via Gau{\ss}-Manin from a smooth,
proper
family of schemes over $X$, we show using work of Katz \cite{Ka} that
$w_n(E,\nabla,P)=0$. A longstanding conjecture of Ogus \cite{Og} would imply
that a
class in $\Gamma(X,\sH^n)$ in characteristic 0, (where  $\sH$ is the Zariski
sheaf
of de Rham cohomology),
which vanished when reduced mod
$p$ for almost all $p$ was $0$. Thus, Ogus' conjecture would imply an
affirmative
answer to \ref{subsubsec:question} for Gau{\ss}-Manin systems.

\subsection{}
In concrete applications, one frequently deals with connections $\nabla$ with
logarithmic poles. Insofar as possible, we develop our constructions in this
context (see section \ref{sec:log}). The most striking remark is that
even if
$\nabla$ has logarithmic poles, $w_n(E, \nabla)$ does not have
any poles (see theorem (\ref{thm:logpoles})).

\subsection{Acknowledgments:}We thank E. Looijenga for making us aware of the
similarity between our invariant
$w_2$ with Witten's invariant for 3-manifolds (see \cite{W}),
and G. van der Geer for
communicating \cite{Ge} to us.

\section{Affine Fibrations}\label{sec:affib}An affine bundle $Y$ over a scheme
$X$ is, by
definition, an $\sV$-torseur for some vector bundle $\sV$. Such things are
classified by
$H^1(X,\sV)$. In particular, Zariski-locally, $Y\cong X\times\A^n$. Pullback
from $X$ to
$Y$ is an isomorphism on Chow motives, and hence on any Weil cohomology; e.g.
$H_{\rm
DR}(X)\cong H_{\rm DR}(Y)$,
$H_{\text{\'et}}(X)\cong H_{\text{\'et}}(Y)$, etc. The following
is known as ``Jouanolou's trick". We recall the argument from \cite{Jo}.
\begin{prop}Let $X$ be a quasi-projective variety. Then there exists an affine
bundle $Y\to
X$ such that $Y$ is an affine variety.
\end{prop}
\begin{proof}Let $X\subset\bar X$ be an open immersion with $\bar X$
projective. Let
$\tilde X$ be the blowup of $\bar X-X$ on $\bar X$. $\tilde X$ is projective,
and
$X\subset\tilde X$ with complement $D$ a Cartier divisor. Suppose we have
constructed $\pi
: \tilde Y\to \tilde X$ an affine bundle with $\tilde Y$ affine. Since the
complement of a
Cartier divisor in an affine variety is affine (the inclusion of the open is
acyclic for
coherent cohomology, so one can use Serre's criterion) it follows that
$\pi^{-1}(X) \to X$
is an affine bundle with $Y:=\pi^{-1}(X)$ affine. We are thus reduced to the
case $X$
projective. Let $P(N) \to \P^N$ be an affine bundle with $P(N)$ affine. Given a
closed
immersion $X\hookrightarrow \P^N$, we may pull back $P(N)$ over $X$, so we are
reduced to
the case $X=\P^N$. In this case, one can take $Y={\rm GL}_{N+1}/{\rm
GL}_{N}\times{\rm
GL}_{N+1}$.
\end{proof}
An exact sequence of vector bundles
$0 \to G \to F\to E \to 0$
on $X$ gives rise to an exact sequence of Hom bundles
$$0 \to Hom(E,G) \to Hom(E,F)\to Hom(E,E) \to 0$$
and so an isomorphism class of affine bundles
$$\partial({\rm Id}_E)\in H^1(X,Hom(E,G)).$$
Of particular interest is the Atiyah sequence. Let $X$ be a smooth variety, and
let
$\sI\subset\sO_X\otimes\sO_X$ be the ideal of the diagonal. Let $\sP_X :=
\sO_X\otimes\sO_X/\sI^2$, and consider the exact sequence
$$0\to \Omega^1_X \to \sP_X\to \sO_X \to 0$$
obtained by identifying $\sI/\sI^2 \cong\Omega^1$ in the usual way. Note that
$\sP_X$ has
two distinct $\sO_X$-module structures, given by multiplication on the left and
right. These
two structures agree on $\Omega^1$ and on $\sO_X$. Given $E$ a vector bundle on
$X$, we
consider the sequence (Atiyah sequence)
\begin{equation} \label{eqn:Atiyah} 0\to E\otimes_{\sO_X}\Omega^1_X \to
E\otimes_{\sO_X}\sP_X\to E
\to 0
\end{equation}
The tensor in the middle is taken using the left $\sO_X$-structure, and then
the sequence is viewed as a sequence of $\sO_X$-modules using the {\it right}
$\sO_X$-structure.
\begin{prop} Connections on $E$ are in $1-1$ correspondence with splittings of
the Atiyah
sequence (\ref{eqn:Atiyah}).
\end{prop}
\begin{proof}(See \cite{At} and  \cite{D}).
As a sequence of sheaves of abelian groups, the Atiyah sequence is split by
$e\mapsto e\otimes 1$. Let $\theta : E\to E\otimes \sP_X$ be an $\sO$-linear
splitting.
Define
$$\nabla(e) := \theta(e) - e\otimes 1 \in  E\otimes\Omega^1_X.$$
We have
\begin{multline*}
\nabla(f\cdot e) := \theta(e)\cdot (1\otimes f) - (e\otimes 1)(f\otimes 1)=\\
(1\otimes f)\cdot \nabla(e) + (e\otimes 1)\cdot (1\otimes f - f\otimes 1) \\
=f\nabla(e)+df\wedge e,
\end{multline*}
which is the connection condition. Conversely, given a connection
$\nabla$, the same
argument shows that
$\theta(e) = \nabla(e)+e\otimes 1$ is an $\sO$-linear splitting.
\end{proof}
\begin{cor}Let $E$ be a vector bundle on a smooth affine variety $X$.
Then $E$ admits an
algebraic connection.
\end{cor}
\begin{proof}An exact sequence of vector bundles on an affine variety admits a
splitting.
\end{proof}
\subsection{}\label{subsec:aff_concl}
In conclusion, given a vector bundle $E$ on a smooth
variety
$X$, there exists two sorts of affine bundles
$\pi : Y\to X$ such that $\pi^*E$ admits a
connection. We can take
$Y$ to be the Atiyah torseur associated to $E$, in which case the connection is
canonical, or
we can take $Y$ to be affine,
in which case all vector bundles admit (non-canonical)
connections.

\section{Chern-Simons}\label{sec:cs} We begin by recalling in an algebraic
context the basic
ideas involving connections and the Chern-Weil and Chern-Simons constructions.
\subsection{Connections and curvature} Let $R$ be a $k$-algebra
of finite type ($R$ and
$k$ commutative with $1$). A connection $\nabla$ on a module
$E$ is a map $\nabla : E \to E\otimes_R\Omega^1_{R/k}$ satisfying
$\nabla(f\cdot e) =
e\otimes df + f\cdot\nabla e$. More generally, if $D\subset\Spec R$ is a
Cartier
divisor, of equation $f$, one defines the module
$\Omega^1_{R/k}(\log D)$
of K\"ahler $1$-forms with logarithmic poles along
$D$,
as the submodule of forms $w$
with poles along $D$ such that $w\cdot f$ and
$w \wedge df$ are regular \cite{DI}.
A connection with log poles along $D$ is a $k$ linear map
$\nabla : E \to E\otimes\Omega^1_{R/k}(\log D)$
fulfilling the Leibniz relations. When $E$ has a global
basis $E=R^N$, $\nabla$ can
be written in the form $d+A$, where $A$ is an $N\times N$-matrix of $1$-forms.
Writing
$e_i = (0,\ldots,1,\ldots,0)$ we have
\[\nabla(e_i) = \sum_je_j\otimes a_{ij}.\]
The map $\nabla$ extends to a map $\nabla : E\otimes\Omega^i \to
E\otimes\Omega^{i+1}$ defined by $\nabla(e\otimes\omega) =
\nabla(e)\wedge\omega + e\otimes d\omega$. The curvature of the connection is
the
map $\nabla^2 : E \to E\otimes\Omega^2$. The curvature is $R$-linear and is
given
in the case $E=R^N$ by
\begin{multline*}
\nabla^2(e_i) = \sum_je_j\otimes da_{ij} + \sum_{j,\ell}e_\ell\otimes
a_{j\ell}\wedge a_{ij} \\= (0,\ldots,1,\ldots,0)\cdot (dA - A^2).
\end{multline*}
The curvature matrix $F(A)$ is defined by $F(A) = dA-A^2$. (Note that the
definition $F(A) = dA+A^2$ is also found in the literature, e.g. in \cite{CS}.)

Given $g\in{\rm GL}_N(R)$, let $\gamma = g^{-1}$. We can rewrite the connection
$\nabla = d+A$ in terms of the basis $\epsilon_i := e_i\cdot g =
(g_{i1},\ldots,g_{iN})$, replacing
$A$ and $F(A)$ by
\begin{eqnarray}
 dg\cdot g^{-1}+gAg^{-1} =& -\gamma^{-1}d\gamma + \gamma^{-1}A\gamma \\
 F(dg\cdot g^{-1}+gAg^{-1}) = &gF(A)g^{-1}.\label{eqn:curv}
\end{eqnarray}
A connection is said to be integrable or flat if $\nabla^2=0$. For a connection
on
$R^N$ this is equivalent to $F(A) = 0$.

\subsection{} We recall some basic ideas from \cite{CS}. Let $\sG$ be a Lie
algebra over a field $k$ of characteristic $0$, and let $G$ be the
corresponding
algebraic group. (The only case we will use is $G={\rm GL}_N$.) Write
$\sG^\ell :=\underbrace{\sG\otimes\cdots\otimes\sG}_{\ell\;\rm factors}$.
$G$ acts
diagonally on $\sG^\ell$ by the adjoint action on each factor, and an element
$P$ in the
linear dual $(\sG^\ell)^*$ is said to be invariant if it is invariant under
this
action. For a $k$-algebra $R$ we consider the
module $\Lambda^{r,\ell} := \sG^\ell\otimes_k\Omega^r_{R/k}$ of $r$-forms on
$R$
with values in $\sG^\ell$. Let $x_i$ denote tangent vector fields, i.e.
elements in
the $R$-dual of $\Omega^1$. We describe two products $\wedge :
\Lambda^{r,\ell}\otimes_R\Lambda^{r',\ell'} \to \Lambda^{r+r',\ell+\ell'}$ and
$[\quad ] : \Lambda^{r,1}\otimes_R\Lambda^{r',1} \to \Lambda^{r+r',1}$.
In terms of
values on tangents, these are given by
\begin{gather}
\varphi\wedge\psi(x_1,\ldots,x_{r+r'}) = \notag \\
\sum_{\pi,{\rm
shuffle}}\sigma(\pi)\varphi(x_{\pi_1},\ldots,x_{\pi_r})
\otimes\psi(x_{\pi_{r+1}},
\ldots,x_{\pi_{r+r'}})\\
{[\varphi , \psi ]}(x_1,\ldots,x_{r+r'}) = \notag \\
\sum_{\pi,{\rm
shuffle}}\sigma(\pi)[\varphi(x_{\pi_1},\ldots,x_{\pi_r}),
\psi(x_{\pi_{r+1}},\ldots,x_{\pi_{r+r'}})]
\end{gather}
Here $\sigma(\pi)$ is the sign of the shuffle. These operations satisfy the
identities (for $P\in(\sG^\ell)^*$ not necessarily invariant)
\begin{eqnarray} [\varphi,\psi] =& (-1)^{rr'+1}[\psi,\varphi] \\
{[[\varphi,\varphi{]},\varphi]} =& 0 \\
d[\varphi,\psi] =& [d\varphi,\psi] + (-1)^r[\varphi,d\psi] \\
d(\varphi\wedge\psi) =& d\varphi\wedge\psi+(-1)^r\varphi\wedge d\psi \\
d(P(\varphi)) =& P(d\varphi) \\
P(\varphi\wedge\psi\wedge\rho) =& (-1)^{rr'}P(\psi\wedge\varphi\wedge\rho)
\end{eqnarray}
If $P$ is invariant, we have in addition for $\varphi_i\in\Lambda^{r_i,1}$ and
$\psi\in\Lambda^{1,1}$
\begin{equation}\sum_{i=1}^\ell
(-1)^{r_1+\ldots+r_i}P(\varphi_1\wedge\cdots\wedge[\varphi_i,\psi]
\wedge\cdots\wedge\varphi_\ell)=0\end{equation}
By way of example, we note that if $A=(a_{ij}), B=(b_{ij})$ are matrices
of 1-forms, then
writing $AB$ (or $A^2$ when $A=B$) for the matrix of $2$-forms with entries
\[ \sum_\ell a_{i\ell}\wedge b_{\ell j} \]
we have
\begin{eqnarray*}
&[A,A](x_1,x_2)_{ij} = \big([A(x_1),A(x_2)] -
[A(x_2),A(x_1)]\big)_{ij} \\
&= 2\big(A(x_1)A(x_2)-A(x_2)A(x_1)\big)_{ij} \\ &= 2\sum_\ell
\big(a_{i\ell}(x_1)a_{\ell j}(x_2) - a_{i\ell}(x_2)a_{\ell j}(x_1)\big) \\
&= 2\sum_\ell a_{i\ell}\wedge a_{\ell j}(x_1,x_2) = 2A^2(x_1,x_2).
\end{eqnarray*}
whence
$$ A^2 = \frac{1}{2}[A,A]
$$
In the following, for $\varphi\in\Lambda^{r,\ell}$ we frequently write
$\varphi^n$
in place of $\varphi\wedge\ldots\wedge\varphi$ ($n$-times). The signs differ
somewhat from \cite{CS} because of our different convention for the curvature
as
explained above.

\begin{thm}{\cite{CS}} Let $P\in(\sG^\ell)^*$ be invariant. To a matrix $A$ of
$1$-forms over a ring $R$, we associate a matrix of $2$-forms depending on a
parameter $t$
$$\varphi_t := tF(A)-\frac{1}{2}(t^2-t)[A,A].$$
Define
\begin{equation}\label{eqn:2.2.10} TP(A) =
\ell\int_0^1P(A\wedge\varphi_t^{\ell-1})dt \in \Omega^{2\ell -1}_{R/k}
\end{equation}
For example, for
\begin{equation}
P(M) = {\rm Tr} M^2, \ \ell = 2, \ TP(A) = {\rm Tr}(AdA -\frac{2}{3}
A^3)
\end{equation}
Then $dTP(A) = P(F(A)^\ell)$.
The association $A\mapsto TP(A)$ is functorial for
maps of rings $R\to S$. If $A\mapsto T'P(A)$ is another such functorial mapping
satisfying $$dT'P(A) = dTP(A) = P(F(A)^\ell),$$ then
$$T'P(A) - TP(A) = d\rho$$ is
exact.
\end{thm}
\begin{proof} The first assertion follows from Prop. 3.2 of \cite{CS}, noting
that $\Omega(A)$ in their notation is $-F(-A)$ in ours. For the second
assertion,
we may assume by functoriality that $R$ is a polynomial ring, so
$H^{2\ell-1}_{DR}(R/k) = (0).$ The form $T'P(A) - TP(A)$ is closed, and hence
exact.
\end{proof}
\begin{prop}\label{prop:basis_ind}
With notation as above, let $g\in {\rm GL}_N(R)$,
and assume
$\ell\ge 2$.  Then $TP(dg\cdot g^{-1}+gAg^{-1}) - TP(A)$
is Zariski-locally exact,
i.e. there exists an open cover $\Spec(R) = \bigcup U_i$ such that the above
expression is exact on each $U_i$.
\end{prop}
\begin{proof} The property of being Zariski-locally exact is compatible under
pullback, so we may argue universally. The matrix $A$ of $1$-forms (resp. the
element $g$) is pulled back from the coordinate ring of some affine space
$\A^m$
(resp. from the universal element in ${\rm GL}_N$ with coefficients in the
coordinate ring of ${\rm GL}_N$),
so we may assume $R$ is the coordinate ring of
$\A^m\times {\rm GL}_N$.

Let $\eta$ be a closed form on a smooth variety $T$. Let
$f : S\to T$ be surjective, with $S$ quasi-projective. Then $\eta$ is
locally exact on $T$ if and only if $f^*\eta$ is locally exact on $S$. Indeed,
given $t\in T$ we can find a section $S''\subset S$ such that the
composition $f' : S' \to T$, where $S' \to S''$
is the normalization, is finite over some neighborhood $t\in U$. Assuming
$f^*\eta$ is locally exact, it follows that $f^*\eta | {f'}^{-1}(U)$ is locally
exact, and so by a trace argument (we are in characteristic
zero) that $\eta | U$ is locally exact as well.

We apply the above argument with
$$\eta = TP(dg\cdot g^{-1}+gAg^{-1}) - TP(A)$$ and
$T = \A^m\times {\rm GL}_N$. As a scheme, ${\rm GL}_N\cong \G_m\times {\rm
SL}_N$,
and for some large integer $M$ we can find a surjection $\coprod_{{\rm
finite}}\A^M
\to {\rm SL}_N$ by taking products of upper and lower triangular matrices with
$1$
on the diagonal and then taking a disjoint sum of translates. Pulling back, it
suffices to show that a closed form of degree $\ge 2$ on $\A^{M+m}\times\G_m$
is
exact. This is clear.
\end{proof}

\subsection{Construction} Let $E$ be a vector bundle of rank $N$ on a smooth
quasi-projective variety $X$. Let $P$ be an invariant polynomial as above of
degree
$n$ on the Lie algebra $\sG\sL_N$. Suppose given a connection $\nabla$ on $E$.
(Such
a connection exists when $X$ is affine because the Atiyah sequence splits) Let
$X =
\bigcup U_i$ be an open affine covering such that $E | U_i \cong \sO^{\oplus
N}$,
and let
$A_i$ be the matrix of $1$-forms corresponding to $\nabla | U_i$. The class of
$TP(A_i)\in \Gamma(U_i,\Omega^{2n-1}/d\Omega^{2n-2})$ is independant of the
choice
of basis for $E | U_i$ by \ref{prop:basis_ind}.
It follows that these classes glue
to give a global class
\begin{equation} w_n(E,\nabla,P) \in \Gamma(X,\Omega^{2n-1}/d\Omega^{2n-2})
\end{equation}

\begin{prop} Let $E$ be a rank $N$-vector bundle on a smooth affine variety
$X$. Let
$\nabla$ and $\nabla'$ be two connections on $E$. Let $P$ be an invariant
polynomial of degree $n$. Then there exists a form
$$\eta\in\Gamma(X,\Omega^{2n-1}_X)$$ such that
$$w_n(E,\nabla,P) - w_n(E,\nabla',P) \equiv \eta\;
{\rm mod} (d\Omega^{2n-2})$$
\end{prop}
\begin{proof} Because $X$ is affine, any affine space bundle $Y\to X$ admits a
section. (An affine space bundle is a torseur under a vector bundle.) Thus, we
may
replace $X$ by an affine space bundle over $X$. Since $X$ is affine, $E$ is
generated by its global sections, so we may find a Grassmannian $G$ and a map
$X\to
G$ such that $E$ is pulled back from $G$. We may find an affine space bundle $Y
\to
G$ with $Y$ affine. Replacing $X$ with $X\times_GY$, which is an affine bundle
over
$X$, we may assume $E$ pulled back from a bundle $F$ on $Y$. Since $Y$ is
affine,
$F$ admits a connection $\Psi$, and it clearly suffices to prove the
proposition
for $\nabla$ the pullback of $\Psi$. Write $\nabla' - \nabla = \gamma$ with
$\gamma\in {\rm Hom}_{\sO_X}(E,E\otimes\Omega^1)$. Let $\iota :
X\hookrightarrow
\A^m$ be a closed immersion. The product map $X\hookrightarrow Y\times
\A^m$ is a
closed immersion, hence $\gamma$ lifts to
$\varphi\in{\rm Hom}_{\sO_{Y\times\A^m}}(F,F\otimes\Omega^1_{Y\times\A^m})$.
Let
$\Psi' := \Psi+\varphi$. We are now reduced to the case $X=Y\times\A^m$.
Writing $\sH^{2n-1}$ for the Zariski cohomology sheaf of the de Rham
complex on $X$,
one knows that $\Gamma(X, \sH^{2n-1})
\subset \Gamma(U, \sH^{2n-1})$ for any
open $U \neq \emptyset$ \cite{BO}. Taking $U = \A^{M+m}$, where
$\A^M$ is an affine cell in  $Y$,
we may assume $\Gamma(X,\sH^{2n-1}) = (0)$. Now $dw(E,\nabla,P)$
and $dw(E,\nabla',P)$ both represent the same class in cohomology, so, since
$X$ is
affine, there exists $\eta\in\Gamma(X,\Omega^{2n-1})$ such that
$$w_n(E,\nabla,P)-w_n(E,\nabla',P)-\eta\in\Gamma(X,\sH^{2n-1}) = (0).$$
\end{proof}

\begin{prop}Let $\nabla$ be an integrable connection on $E$, and let $P$ be an
invariant
polynomial of degree $n$. Let $\sH^{2n-1} = \Omega^{2n-1}_{\rm
closed}/d\Omega^{2n-2}$. Then
$w_n(E,\nabla,P)\in\Gamma(X,\sH^{2n-1})$, i.e. $dw=0$.
\end{prop}
\begin{proof}$dw = P(F(\nabla))=0$ since $\nabla$ integrable implies
$F(\nabla)=0$.
\end{proof}

\begin{prop}\label{prop:unicity}Let $\nabla$ be an integrable connection on
$E$, and
let $P = \lambda P_n + Q$ be an invariant polynomial of degree
$n$, where $P_n$ is the $n$-th elementary symmetric function
and
$$Q= \sum \mu_{ij}P_i\cdot P_j, i\neq j, i\geq 1, j\geq 1, \lambda \in \Q,
\mu_{ij} \in \Q.
$$
Then $$w_n(E, \nabla, P) =\lambda w_n(E,\nabla)$$ (see notation
(0.2.3)).
\end{prop}
\begin{proof}It is enough to see that $w_n(E, \nabla, P_i \cdot P_j) =0$
for $i\geq 1, j\geq 1$. One has $P_j(F(A)^j) =0$, so by
Theorem 2.2.1 $$T'(P_i \cdot P_j)(A) = TP_i(A) \cdot P_j((F(A_j))^j)=0$$
differs from $T(P_i \cdot P_j)(A)$ by an exact form on the open on
which $A$ is defined.
\end{proof}

\subsection{Rigidity}
\begin{thm}\label{prop:rigid}
Let $f: X \to S$ be a smooth proper morphism between smooth algebraic
varieties defined over a field $k$ of characteristic zero. Assume ${\rm dim }
S = 1$.
Let $\nabla : E \to \Omega^1_{X/S} \otimes E$ be a relative flat connection,
and $P$ be an invariant polynomial. Then $$w_n(E, \nabla, P) \in
H^0(X, \sH^{2n-1}(X/S))$$ lifts canonically to a class in $H^0(X, \sH^{2n-1})$
for
$n \geq 2$.
\end{thm}
\begin{proof}
Take locally the matrix $A'_i \in H^0(X_i, M(N, \Omega^1_{X/S}))$ of the
connection, $N$ being the rank of $E$. Take liftings
$A_i \in H^0(X_i, M(N, \Omega^1_{X}))$, and define
$TP(A_i)$ looking at the $\Omega^1_X$ valued connection defined by
$A_i$. Since $F(A_i) \in H^0(X_i, M(N, f^*\Omega^1_S))$, one
has $F(A_i)^n = 0 \ {\rm for} \ n \geq 2$, and $dTP(A_i) =
P(F(A)^n) = 0$. On $X_i \cap X_j$, one has
$$A_j = dg \cdot g^{-1} + g A_i g^{-1} - \Gamma_{ij}$$
where $\Gamma_{ij} \in H^0(X_i \cap X_j, f^*\Omega^1_S)$.
Using proposition 2.2.2, we just have to show
that $TP(B) - TP(B + \Gamma)$ is locally exact
for some matrix of one forms $B =
dg \cdot g^{-1} + g  A_i  g^{-1}$,
verifying $F(B)  \omega = 0$ for any $w \in  M(N, f^*\Omega^1_S)$,
and $\Gamma = \Gamma_{ij} \in f^* \Omega^1_S$.
By \ref{prop:unicity} it is enough to consider $P(M) = {\rm Tr}M^n$.
One has
\begin{equation}\varphi_t(B + \Gamma) = F(t(B + \Gamma)) =
F(tB) + td \Gamma -t^2(\Gamma B + B \Gamma)
\end{equation}
and
\begin{equation}{}
F(tB) \omega = (t-t^2) dB \omega
\end{equation}
with $\omega$ as above. Thus
\begin{gather}{}
P((B+\Gamma) \wedge \varphi_t^{n-1}(B + \Gamma)) = \\
{\rm Tr} (B+\Gamma)[(tdB-t^2B^2)^{n-1} + (n-1)
(t-t^2)^{n-2}(dB)^{n-2} (t d \Gamma -t^2(B\Gamma + \Gamma B)] \notag \\
= P(B \wedge \varphi_t^{n-1}(B)) + R \notag
\end{gather}
with
\begin{gather}{}
R = {\rm Tr} \Gamma (dB)^{n-1} [ (t-t^2)^{n-1} -2 t^2
(n-1)(t-t^2)^{n-2}] \notag \\
+ (n-1)(t-t^2)^{n-2}t {\rm Tr}B(dB)^{n-2} d\Gamma
\end{gather}
Write
${\rm Tr} d(B \Gamma) = {\rm Tr} dB \Gamma - {\rm Tr} Bd
\Gamma$. Then we have
\begin{equation}{}
R = F(t) {\rm Tr} \Gamma (dB)^{n-1}
\end{equation}
modulo exact forms, with
\begin{equation}{}
F(t) = n(t-t^2)^{n-1} - (n-1)t^2(t-t^2)^{n-2} =
(t(t-t^2)^{n-1})'.
\end{equation}
The assertion now follows from \eqref{eqn:2.2.10}.
\end{proof}

\section{Flat Bundles}\label{sec:flat}
The following notations will reoccur frequently.
\subsection{}$X$ will be a
smooth variety, and
$D=\bigcup D_i\subset X$ will be a normal crossings divisor, with $j : X-D \to
X$. We
will assume unless otherwise specified that the ground field $k$ has
characteristic
$0$.
\subsection{}\label{subsec:bundle}$(E,\nabla)$ will be a vector bundle $E$ of
rank
$r$ on
$X$ with connection $\nabla : E \to E\otimes\Omega^1_{X/k}({\rm log}\;D)$
having
logarithmic poles along $D$. The Poincar\'e residue map $\Omega^1({\rm log}\;D)
\to
\sO_{D_i}$ is denoted ${\rm res}_{D_i}$, and $\Gamma_i := {\rm
res}_{D_i}\circ\nabla
: E \to E|_{D_i}$.
\subsection{}When $E$ is trivialized on the open cover $X = \cup
X_i$, with basis $\underline{e}_i$ on $X_i$, then $(E, \nabla)$
is equivalent to the data
\begin{eqnarray*}
g_{ij} &\in &\Gamma (X_i \cap X_j, GL (r, \sO_X)) \\
g_{ik} &=& g_{ij} g_{jk} \\
A_i &\in &\Gamma (X_i, M (r, \Omega^{1}_{X} ({\rm log} D)))
\end{eqnarray*}
with $g^{-1}_{ij} dg_{ij} = g^{-1}_{ij} A_i g_{ij} - A_j$.
\subsection{}The curvature
$$
\nabla^2 : E \to \Omega^{2}_{X} ({\rm log} D) \otimes E
$$
is given locally by
$$
\nabla^2 = F (A_i) : = d A_i - A_i A_i.
$$
The connection $\nabla$ is said to be flat, or integrable if $\nabla^2=0$.

\subsection{}For two $r \times r$ matrices $A$ and $B$ of differential forms of
weight $a$ and $b$ respectively, one writes ${\rm Tr} \ AB$ for the trace of
the
$r \times r$ matrix $AB$ of weight $a+b$, and one has ${\rm Tr} \
AB = (-1)^{ab} {\rm Tr} \ BA$. We denote by ${}^t A$ the transpose of $A$:
$({}^t A)_{ij} = A_{ji}$.

\subsection{} For any cohomology theory $H$ with a localization
sequence, the $i$-th level of Grothendieck's coniveau
filtration is defined by
\begin{multline*}N^i H^{\bullet} = \{ x \in H^{\bullet} | \exists
\mbox{  subvariety } Z \subset X \mbox { of codimension } \geq i \\
\mbox{ such that }
0 = x |_{X-Z} \in H^{\bullet} (X-Z).
\end{multline*}

\subsection{}For any cohomology theory $H$ defined in a topology
finer than the Zariski topology, one defines the Zariski sheaves
$\sH$ associated to the presheaves $U\mapsto H(U)$ (\cite{BO}). When $H$ is the
cohomology for the analytic topology with coefficients in a constant sheaf $A$,
we
sometimes write $\sH(A)$.
 For example the
Betti or de Rham sheaves $\sH(\C)$ are simply the cohomology sheaves for the
complex
of algebraic differentials $\Omega^*_X$. For $D\subset X$ as above, we write
$\sH^{\bullet}(\log D)$ for the cohomology sheaves of $\Omega^*(\log D)$. It is
known that $\sH^{\bullet}(\log D)\cong j_*\sH^\cdot_{X-D}$.

\subsection{}When $k=\C$ we use the same notation $\Omega^*_X$ for the analytic
and algebraic de Rham complexes. For integers $a$ and $b$, the analytic Deligne
cohomology is defined to be the hypercohomology of the complex of analytic
sheaves
$$H^a_{\sD,{\rm an}}(X,\Z(b)) := \H^a(X_{{\rm an}}, \Z(b) \to \sO \to
\Omega^1\to\cdots\to \Omega^{b-1}).
$$
(This should be distinguished from the usual Deligne cohomology, which is
defined using differentials with at worst log poles at infinity.) One has a
cycle
class map from the Chow group of algebraic cycles modulo rational equivalence
to Deligne cohomology:
$$CH^i(X) \to H^{2i}_{\sD,{\rm an}}(X,\Z(i)).$$

\subsection{}We continue to assume $k=\C$. Let $\alpha : X_{{\rm an}} \to
X_{{\rm
Zar}}$ be the identity map. For a complex $C^\cdot$, let $t_{\ge i}C^\cdot$ be
the
subcomplex which is zero in degrees $<i$ and coincides with $C^\cdot$ in
degrees
$\ge i$. There is a map of complexes $t_{\ge i}C^\cdot \to C^\cdot$. The
complex
$$\Z(j) \to \sO_X \to \Omega^1_X\to \cdots$$
in the analytic topology is quasi-isomorphic to the cone $\Z(j)\to \C$, and
hence
to $\C/\Z(j)[-1]$. We obtain in this way a map in the derived category
$$(t_{\ge j}\Omega^\cdot_X) \to \C/\Z(j)$$
The kernel of the resulting map
$$R^{j}\alpha_*(t_{\ge j}\Omega^\cdot_X) \cong \ker\big(\alpha_*\Omega^j\to
\alpha_*\Omega^{j+1}\big) \to R^{j}\alpha_*(\C/\Z(j))
 $$
is denoted $\Omega^j_{\Z(j)}$ \cite{EII}.
Note $\Omega^j_{\Z(j)}$ is a Zariski sheaf. Writing $\sK^m_j$ for the Milnor
$K$-sheaf (subsheaf of the constant sheaf $K^{{\rm Milnor}}_j(k(X))$), the
$d\log$-map
$$\{f_1,\ldots,f_j\} \mapsto df_1/f_1\wedge\cdots\wedge df_j/f_j$$
induces a map
\begin{equation}\label{eqn:K} d\log : \sK^m_j \to
\Omega^j_{\Z(j)}.\end{equation}
To see this, note the exponential sequence induces a map
$$\sO_{X_{{\rm Zar}}}^* \to R^1\alpha_*\Z(1)
$$
and we get by cup product a commutative diagram with left hand vertical arrow
surjective
\begin{equation*}
\begin{CD}
\sO_{X_{{\rm Zar}}}^{*\otimes j} \>>> R^j\alpha_*\Z(j) \\
\V V {\rm surj.} V                                     \V V V \\
\sK^m_j \>d\log >>             R^j\alpha_*\C.
\end{CD}
\end{equation*}

We shall need some more precise results about the sheaf $\Omega^j_{\Z(j)}$.
\begin{lem}\label{lem:maps}
\begin{enumerate}
\item \label{first}There is a natural map
$$H^i(X_{{\rm Zar}},\Omega^i_{\Z(i)})\to H^{2i}_{\sD,{\rm an}}(X,\Z(i))
$$
\item \label{second}Let $D\subset X$ be a normal crossings divisor. Then there
is a
natural map
\begin{multline*}\H^i(X_{{\rm
Zar}},\Omega^i_{\Z(i)}\to\alpha_*\Omega^i_X(\log\;D)\to\ldots)\to \\
\H^{2i}(X_{\rm
an},\Z(i)\to {\sO}_X \to \Omega^1_X\to\ldots\to\Omega^{i-1}_X
\to\Omega^i(\log\;D)_X\to\ldots)
\end{multline*}
\item There is a natural map
\begin{multline*}
\varphi:\H^i(X_{{\rm
Zar}},\sK_i^m \xrightarrow{d\log}\Omega^i_X(\log\;D)\to\ldots)\to \\
\H^{2i}(X_{\rm
an},\Z(i)\to {\sO}_X \to \Omega^1_X\to\ldots\to\Omega^{i-1}_X
\to\Omega^i(\log\;D)_X\to\ldots)
\end{multline*}
In particular, for $D=\emptyset$, we get a map
$$\H^i(X_{{\rm Zar}}, \sK_i^m \xrightarrow{d\log} \Omega^i_X\to\ldots)\to
H^{2i-1}(X_{\rm an},\C/\Z(i)).$$
\end{enumerate}
\end{lem}
\begin{proof}We consider the spectral sequence
\begin{gather*}
R^j :=
R^j\alpha_*(\Z(i)\to\sO\to\Omega^1\to\ldots\to\Omega^{i-1}) \\
E_2^{p,q} = H^p(X_{{\rm Zar}},
R^q)\Rightarrow H^{p+q}_{\sD,{\rm an}}(X,\Z(i))
\end{gather*}
One checks that
\begin{gather*}
R^s \cong \sH^{s-1}(\C/\Z(i));\qquad s< i \\
0 \to \sH^{i-1}(\C/\Z(i)) \to R^i \to \Omega^i_{\Z(i)} \to 0 \\
0 \to \sH^{i-1}(\C/\Z(i)) \to R^s \to \ker(\sH^s(\C)\to\sH^s(\C/\Z(i))) \to
0;\qquad s>i
\end{gather*}
We have by (\cite{BO}) that $H^a(X_{\rm Zar},\sH^b(A))=(0)$ for $a>b$ and $A$
any
constant sheaf of abelian groups. Applying this to the above, we conclude
$E_2^{a,2i-a}=H^a(X_{\rm Zar},R^{2i-a})=(0)$ for $a>i$, and $E_2^{i,i}\cong
H^i(X,\Omega^i_{\Z(i)})$. Assertion (\ref{first}) follows. The construction of
the
map in (\ref{second}) is similar and is left for the reader. Finally, (3)
follows
by composing the arrow from (2) with the $d$ log map \ref{eqn:K}.

\end{proof}
\subsection{Characteristic classes.}Let $(E,\nabla)$ be a bundle with
connection as
in \ref{subsec:bundle} and assume $\nabla$ is flat. Functorial and additive
characteristic classes
$$c_i(E,\nabla)\in \H^i(X_{{\rm Zar}},\sK^m_i\to
\Omega^i_X(\log D)\to\Omega^{i+1}_X(\log D)\to\ldots)
$$
were defined in \cite{EI}. These classes have the following compatibilities:
\subsubsection{} Under the map
\begin{multline*}\H^i(X_{{\rm Zar}},\sK^m_i\to
\Omega^i_X(\log\;D)\to\Omega^{i+1}_X(\log\;D)\to\ldots) \to \\
H^i(X,\sK^m_i)\cong
CH^i(X)
\end{multline*}
we have $c_i(E,\nabla)\mapsto c_i^{{\rm Chow}}(E)\in CH^i(X)$.
\subsubsection{}Assume $X$ proper and $D=\phi$. The classes $c_i(E,\nabla)$
lift classes
$c_i^{\rm an}(E,\nabla)\in H^{2i-1}(X_{\rm an},\C/\Z(i))$ defined in
\cite{EII}, via the
commutative diagram
\begin{equation*}\begin{CD}
\H^i (X, \sK^{m}_{i} \to \Omega^{i}_{X} \to \Omega^{i+1}_{X} \to
\ldots) \>>> CH^i (X) \\
\V \varphi\;(\ref{lem:maps}(3)) V V \V V \psi={\rm cycle\; map} V \\
H^{2i-1} (X_{{\rm an}}, \C /\Z (i)) \>>> H^{2i}_{\sD} (X, i).
\end{CD}\end{equation*}
\subsubsection{}When $D\ne \phi$ and $X$ is proper, classes
$$c_i^{\rm an}(E,\nabla)\in \H^{2i} (X_{\rm an}, \Z (i)
\to \sO_X \to \ldots\to \Omega^{i-1}_{X} \to
\Omega^{i}_{X} ({\rm log} D) \to \ldots)
$$
lifting $c_i^{\sD}(E)\in H^{2i}_\sD(X,\Z(i))$ are defined in \cite{EII}. In
general, for
$X$ not proper, these classes lift
$$c^{\sD}_{i} (E|_{X -D} ) \in H^{2i}_{\sD} (X-D, \Z (i))$$
via the factorization through $H^{2i-1} (X-D, \C/\Z (i))$ (\cite{ECon}, (3.5)).
\begin{prop}\label{prop:3.10.1}
The map $\varphi$ from (\ref{lem:maps}(3)) carries $c_i(E,\nabla)$ to
$c_i^{\rm an}(E,\nabla)$. For $X$ proper, the diagram
\begin{equation*}\minCDarrowwidth5pt\begin{CD}
\H^i (X, \sK^m \to \Omega^{i}_{X} ({\rm log} D) \to
\Omega^{i+1}_{X} ({\rm log} D) \to \ldots ) \>>>  CH^i (X) \\
\V \varphi V V  \V V \psi V \\
\H^{2i} (X_{{\rm an}}, \Z (i) \to \sO_X \to \ldots \Omega^{i-1}_{X}
\to\Omega^{i}_{X} ({\rm log} D)\to \ldots ) \>>> H^{2i}_{\sD} (X,
\Z (i))
\end{CD}\end{equation*}
commutes. For $X$ not proper, the diagram remains commutative if one replaces
the bottom
row by
$$
H^{2i-1} ((X -D)_{\rm an}, \C /\Z (i)) \to H^{2i}_{\sD} (X - D, \Z(i))
$$
or if one replaces $H^{2i}_{\sD} (X, \Z (i))$ by $H^{2i}_{\sD,{\rm an}} (X, \Z
(i))$.
\end{prop}
\begin{proof}The central point, for which we refer the reader to (\cite{EII})
is the
following. Let
$\pi : G\to X$ be the flag bundle of $E$ over which $E$ has a filtration
$E_{i-1} \subset
E_i$ by
$\tau \nabla$ stable subbundles with successive rank 1 quotients $(L_i, \tau
\nabla)$
(see \cite{EI}). Then $c_i (E, \nabla)$ and $c_i^{\rm an} (E, \nabla)$
are both defined on $G$ by products starting from
$$
c_1 (L_{\alpha} , \tau \nabla) \in \H^1 (G, \sK_1 \to \pi^*
\Omega^{1}_{X} ({\rm log} D) \to \ldots )
$$
$$
c_1^{\rm an} (L_{\alpha}, \tau \nabla) \in \H^2 (G, \Z (i) \to \sO_G \to
\pi^* \Omega^{1}_{X} ({\rm log} D) \to \ldots ).
$$
It suffices to observe that the ``algebraic'' product
\begin{multline*}
\H^1 (G, \sK_{1} \to \pi^* \Omega^{1}_{X} ({\rm log} D)
 \to \ldots )^{\otimes i} \\
 \to \H^i (G, \sK^{m}_{i} \to \pi^*
\Omega^{i}_{X} (\log D) \to \ldots )
\end{multline*}
(\cite{EII}, p. 51) is defined compatibily with the ``analytic''
product
\begin{multline*}
\H^2 (G, \Z (i) \to \sO_G \to \pi^* \Omega^{1}_{X} ({\rm log} D)
\to \ldots )^{\otimes i} \\
\to \H^{2i} (G, \Z (i) \to \sO_G \to \ldots \to \Omega^{i-1}_{G}
\to \pi^* \Omega^{i}_{X} ({\rm log} D) \to \ldots ).
\end{multline*}
\end{proof}

\subsection{}\label{subsec:N}
Let $\tau:\Omega^*_X \to N^*$ be a map of complexes , with $\sO_X = N^0$,
such that
if $a$ is the smallest degree $b$ for which $B^b :=
{\rm Ker} \Omega^b_X \to N^b \neq 0$, then $B^b = B^a \wedge \Omega^{b-a}_X$.
For example, let $\nabla : \sF \to \Omega^1_X( \log D)\otimes \sF$ be a non
integrable
connection. Then the local relation $dF(A) = [A, F(A)]$ shows that one
can define $N^*$ by
\begin{gather}{}
N^1 = \Omega^1_X( \log D) \\
N^i = \Omega^i_X( \log D)/B^2 \wedge \Omega^{i-2}_X(\log D) \notag \\ \notag
\end{gather}
where $B^2$ is locally generated by the entries of the curvature
matrix of $\nabla$.

Let $(E, \nabla)$ be a flat $N^*$ valued
connection, that is a $k$ linear map
$\nabla : E \to N^1 \otimes E$ satisfying
the Leibniz rule
\begin{equation}
\nabla(\lambda e) = \tau d \lambda e + \lambda \nabla(e),
\end{equation}
the sign convention
\begin{equation}
\nabla(\omega \otimes e) = \tau d (\omega) \otimes e + (-1)^o \omega
\wedge \nabla(e),
\end{equation}
where $o = {\rm deg} \omega$, and $(\nabla)^2 = 0$.
Then the computations of \cite{EI} and \cite{EII}
allow to show the existence of functorial and additive classes
\begin{equation}
c_i(E, \nabla) \in \H^i(X, \sK_i^m \to N^i \to N^{i+1} ...)
\end{equation}
mapping to analytic classes
\begin{equation}
c_i^{{\rm an}}(E, \nabla) \in \H^{2i}(X_{{\rm an}}, \Z(i) \to ...\to
\Omega^i_X \to N^{i+1} ..)
\end{equation}
compatibly with the classes $c_i^{\sD}(E)$ and $c_i^{{\rm Chow}}(E)$ as before.
As we won't need those classes, we don't repeat in details the construction.

\subsection{}Finally, the $c_i(E,\nabla)$ map to classes
$$
\theta_i (E, \nabla) \in H^0 (X, \frac{\Omega^{2i-1}_{X}
({\rm log} D)}{d \Omega^{2i-2}_{X} ({\rm log} D)}).
$$
In the next section these will be related to the classes $w_i(E,\nabla)$.

\section{The classes $\theta_n$ and $w_n$}\label{sec:thetaw}

Recall that we had defined $w_n(E, \nabla) = w_n(E, \nabla, P_n)$
in (0.2.3) for the $n$-th elementary symmetric function $P_n$.

\begin{thm}\label{thm4.0.1} Let $X$ be a smooth quasi-projective variety over
$\mathbb C$. Let
$E,\nabla$ be a rank
$d$ vector bundle on $X$ with integrable connection. For $d\ge n\ge 2$,
we have $w_n(E,\nabla) =
\theta_n(E,\nabla)$, and $w_n(E, \nabla, P) = \lambda \theta_n(E, \nabla)$
(with the notations of Proposition 2.3.2).
\end{thm}

The proof will take up this entire section. We begin with

\begin{remark} We may assume $X$ is affine, and
$E\cong \mathcal O_X^{\oplus N}$. In this situation,
the class $w_n(E,\nabla)$ lifts canonically to a class in \linebreak
$H^0(X,\Omega^{2n-1})/dH^0(X,\Omega^{2n-2})$.
\end{remark}

Indeed, one knows from \cite{BO} that for $U\subset X$ non-empty open, the
restriction map
$H^0(X,\mc H^{2n-1})\to H^0(U,\mc H^{2n-1})$ is injective.
The assertion about lifting follows from
the construction of $w_n$ in section \ref{sec:cs}
because the trivialization can be taken globally.

\subsection{} The connection is now given by a matrix of $1$-forms and so can
be pulled back in many
ways from some (non-integrable) connection $\Psi$ on the trivial bundle
$\mathcal E\cong \mathcal
O_{\Bbb A^p}^N$. We will want to assume $\Psi$ ``general" in a sense to be
specified below. For
convenience, write $T=\Bbb A^p$ and let $\varphi : X\to T$ be the map pulling
back the
connection. Let $\pi : P \to T$ be the flag bundle for $\mathcal E$ and let
$Q=\varphi^*P$, so we
get a diagram
\begin{equation}\label{eqn4.1.1}\begin{CD} Q \>\varphi >> P \\
\V V\pi V    \V V \pi V \\
X \>\varphi >> T.
\end{CD}
\end{equation}

\subsection{} The curvature $F(\Psi)$ defines an $\mc O_T$-linear map
\begin{equation} F(\Psi) : \mc E \to \mc E\otimes_{\mc O_T} \Omega^2_T.
\end{equation}
In concrete terms, we take $p=2N^2q$ for some large integer $q$, and we write
$x^{(k)}_{ij}$ and
$y^{(k)}_{ij}$ for $1\le i,j\le N$ and $1\le k\le q$ for the coordinates on
$\mathbb A^p$. The
connection $\Psi$ then corresponds to an $N\times N$ matrix of $1$-forms
$A=(a_{ij})$, and the
curvature is given by $F(\Psi) :=(f_{ij})= dA-A^2$. We take
\begin{equation}  a_{ij} =
\sum_{\ell=1}^q x^{(\ell)}_{ij}dy^{(\ell)}_{ij};\quad
f_{ij}=da_{ij}-\sum_{m=1}^Na_{im}\wedge a_{mj}.
\end{equation}
Notice that for $q$ large, we can find $\varphi : X\to T$ so that
$(\mc E,\Psi)$ pulls back to
$(E,\nabla)$.

\subsection{}\label{subsec:M}
We want to argue universally by computing
characteristic classes for $(\mc E,\Psi)$,
but the curvature gets in the way. We could try to kill the
curvature and look for classes in the
quotient complex of $\Omega^*_T$
modulo the differential ideal generated by the $f_{ij}$
(see \ref{subsec:N}), but this
gratuitous violence seems to lead to difficulties. Instead, we will use the
notion of
$\tau$-connection defined in \cite{EI} and \cite{EII} and work with a sheaf of
differential algebras
\begin{equation} M^* = \Omega^*_P/\mc I
\end{equation}
on the flag bundle $P$.

Let
\begin{equation}\label{eqn4.3.2} \mu : \Omega^1_{P/T}\xrightarrow{\iota}
\pi^*Hom(\mc E,\mc E)
\xrightarrow{\pi^*F(\Psi)} \pi^*\Omega^2_T
\end{equation}
be the composition, where $\iota$ is the standard inclusion on a flag bundle.
An easy way to see $\iota$ is to consider the
fibration $ R \to P = R/B$, where $R$ is the corresponding
principal $G = GL(N)$ bundle and $B$ is the Borel
subgroup of upper triangular matrices, and to write
the surjection $\sT(R/T)/B \to \sT(P/T)$ dual to $\iota$,
where $\sT(A/B)$ is the relative tangent space of $A$ with
respect to $B$.
There is an induced map
of graded $\pi^*\Omega^*_T$-modules,
and we define the graded algebra $M^*$ to be the cokernel as
indicated:
\begin{equation}\label{eqn4.3.3}
\Omega^1_{P/T}\otimes_{\mc O_P}\pi^*\Omega^*_T[-2]
\xrightarrow{\mu\otimes 1}
\pi^*\Omega^*_T\to M^*
\to 0.
\end{equation}
Note $M^0=\mc O_P$ and $M^1=\pi^*\Omega^1_T$.

\refstepcounter{subsection}
\begin{prop}\label{prop4.4.1} \begin{enumerate}\item[(i)]
Associated to the connection $\Psi$ on $\mc
E$ there is an
$\mc O_P$-linear splitting $\tau : \Omega^1_P \twoheadrightarrow
\pi^*\Omega^1_T$ of the natural
inclusion $\pi^*\Omega^1_T \xrightarrow{i} \Omega^1_P$. The resulting map
$\delta := \tau\circ d :
\mc O_P \to \pi^*\Omega^1_T$ is a derivation, which coincides with the exterior
derivative on
$\pi^{-1}\mc O_T\subset\mc O_P$. By extension, one defines
$$\delta : \pi^*\Omega^n_T \to \pi^*\Omega^{n+1}_T;\quad
\delta(f\pi^{-1}\omega) =
f\pi^{-1}d\omega+\delta(f)\wedge \pi^{-1}\omega
$$
\item[(ii)] One has
$$\delta^2 = \mu\circ d_{P/T} : \mc O_P \xrightarrow{d_{P/T}}
 \Omega^1_{P/T} \to \pi^*\Omega^2_T,
$$
where $\mu$ is as in \eqref{eqn4.3.2}.
\item[(iii)] There is an induced map $\delta : M^n \to M^{n+1}$ making $M^*$ a
differential graded
algebra. The quotient map $\Omega^*_P \twoheadrightarrow \pi^*\Omega^*_T
\twoheadrightarrow M^*$ is
a map of differential graded $\mc O_P$-algebras.
\end{enumerate}
\end{prop}

\begin{proof} We will give a somewhat different construction of $M^*$ which we
will show coincides
with that defined by \eqref{eqn4.3.3}.

Let $Y$ be a scheme, and let $\mc F$ be a vector bundle on $Y$. Let $\pi_1 :
P_1:=
\mb P(\mc F)\to Y$. Let
$\mc I\subset\Omega^*_{Y}$ be a differential graded ideal, and write $M^*_0 =
\Omega^*_{Y}/\mc
I$. (All our differential graded ideals will be trivial in degree $0$, so
$M^0_0 = \mc O_{Y}$.)
Assume we are given an $M_0$-connection
$\nabla: \mc F \to \mc F\otimes M^1_0$ in an obvious
sense, that is a $k$ linear map fulfilling the ''Leibniz'' rule
$\nabla(\lambda f) = \delta (\lambda) f + \lambda \nabla (f) $
for $\lambda \in \sO_Y, f \in \sF$.
Define $\mc J := \pi_1^{-1}\mc I\Omega^*_{P_1}\subset\Omega^*_{P_1}$, and let
$\tilde M^* :=
\Omega^*_{P_1}/\mc J$. As a consequence of the Leibniz rule,
the pullback $\pi_1^*\mc F$ has a $\tilde M^*$-connection $\tilde\nabla
: \pi_1^*\mc F \to \pi_1^*\mc F \otimes \tilde M^1$.

We want to construct a quotient differential graded algebra $\tilde M^*
\twoheadrightarrow M_1^*$
such that with respect to the quotient $M_1$-connection, the universal sequence
\begin{equation}\label{eqn4.4.1} 0 \to \Omega^1_{P_1/Y}(1) \xrightarrow{j}
\pi_1^*\mc F
\xrightarrow{q}
\mc O_{P_1}(1)
\to 0
\end{equation}
is horizontal. The composition
\begin{equation}\label{eqn4.4.2}\Omega^1_{P_1/Y}(1) \xrightarrow{j} \pi_1^*\mc
F\xrightarrow{\tilde\nabla}\pi_1^*\mc F\otimes\tilde M^1\xrightarrow{q\otimes
1} \tilde M^1(1)
\end{equation}
is easily checked to be $\mc {O}_{P_1}$-linear.
Let $\tilde{\mc {K}_1}\subset \tilde M^1$ denote the image
of the above map twisted by $\mc {O}_{P_1}(-1)$. Define
$\tilde {\mc {K}^*}\subset\tilde M^*$ to be the
graded ideal generated by $\tilde{\mc {K}^1}$ in degree $1$ and
$\delta\tilde{\mc {K}^1}$ in degree $2$. Let
$M_1^* := \tilde{\mc {M}^*}/\tilde{\mc {K}^*}$. It is immediate that $M_1^*$ is
a differential graded
algebra, and that the subbundle $\Omega^1_{P_1/Y}(1)\subset\pi_1^*\mc F$ is
horizontal for the
quotient connection $\pi_1^*\mc F\to \pi_1^*\mc F\otimes M_1^1$.

Now let $P$ denote the flag bundle for $\mc F$. Realize $P$
as a tower of projective bundles
$$ P=P_{N-1}\to\ldots\to P_2\to P_1\to Y
$$
where $P_i$ is the projective bundle on the tautological subbundle on
$P_{i-1}$. Starting with an
$M_0$-connection on $\mc F$ on $Y$, we can iterate the above construction to
get a sheaf of
differential graded algebras $M_i^*$ on $P_i$, and an $M_i$-connection on $\mc
F|P_i$ such that the
tautological partial flag is horizontal. Let $M'{}^*$ be the resulting sheaf of
differential graded
algebras on $P$.

Suppose $M_0^1 = \Omega^1_Y$. We will show by induction on $i$ that $M_i^1\cong
\Omega^1_Y|P_i$ in
such a way that the surjection $\Omega^1_{P_i} \twoheadrightarrow M_i^1$ splits
the natural
inclusion $\Omega^1_Y|P_i\hookrightarrow \Omega^1_{P_i}$, or in other words
that the kernel of the
former is complementary to the image of the latter. This assertion is local on
$P_{i-1}$ (in fact,
it is local on $P_i$), so we may assume $P_i=\mb P(\mc G)$ where $\mc G$ is
trivial on $P_{i-1}$. We
can then lift the $M_{i-1}$-connection on $\mc G$ to an
$\Omega^1_{P_{i-1}}$-connection. The analog of
\eqref{eqn4.4.2} is now
$$\Omega^1_{P_i/P_{i-1}}(1) \to \mc G|P_i \to \mc G|P_i\otimes
\Omega^1_{P_i}\to\Omega^1_{P_i}(1)
$$
This composition twisted by $\mc O_{P_i}(-1)$ is shown in \cite{EI} (0.6.1) to
be (upto sign) a
splitting of
$\Omega^1_{P_i}\to \Omega^1_{P_i/P_{i-1}}$. In particular, its image is
complementary to
$\Omega^1_{P_{i-1}}|P_i$. Factoring out $\Omega^1_{P_i}$ by the image of this
map and by the
pullback of the kernel of $\Omega^1_{P_{i-1}}\twoheadrightarrow M^1_{i-1}\cong
\Omega^1_Y|P_{i-1}$,
it follows easily that $M^1_i\cong \Omega^1_Y|P_i$ as claimed.

To show $M'{}^*$ as constructed here coincides with $M^*$ from \eqref{eqn4.3.3}
we must prove for
$Y=T$ and $\mc F=\mc E$ that $M^2\cong M'{}^2$. We filter $\Omega^1_{P/T}$ so
$\text{fil}_0 =
(0)$ and $\text{gr}_i =\Omega^1_{P_i/P_{i-1}}|P$. We will show by induction on
$i$ that with
reference to \eqref{eqn4.3.3} we have
\begin{equation}\label{eqn4.4.3}
\mu(\text{gr}_i\Omega^1_{P/T})=\delta(\mc
K_{i}^1)\subset(\Omega^2_T|P)/\mu(\text{fil}_{i-1}\Omega^1_{P/T})
\end{equation}
where $\mc K_i^1$ is the image of $\Omega^1_{P_i/P_{i-1}}|P$ in $\Omega^1_T|P$
under the map
analogous to \eqref{eqn4.4.2}. Suppose first $i=1$. Let $e_0, e_1,\ldots$ be a
basis of $\mc E$, and
let
$t_i$ be the corresponding homogeneous coordinates on $P_1 = \mb P(\mc E)$ so
$q(e_i)=t_i$ in
\eqref{eqn4.4.1}. The inclusion $j : \Omega^1_{P_1/T}(1)\hookrightarrow
\pi_1^*\mc E$ is given by
\begin{equation}\label{eqn4.4.4} t_0d(t_i/t_0) \mapsto e_i - (t_i/t_0)e_0.
\end{equation}
Consider the diagram
\begin{equation}\label{eqn4.4.5}
\begin{array}{ccrcc}
\Omega^1_{P_1/T}(1) & \twoheadrightarrow & \mc K^1(1) \subset \Omega^1_{P_1}(1)
&\stackrel{a}{\dashrightarrow} & (\Omega^2_{P_1}/\mc K^1\wedge\Omega^1)(1) \\
j\downarrow &&\uparrow\qquad    &&\uparrow q\otimes 1\\
\pi_1^*\mc E & \xrightarrow{\pi_1^*\Psi} & \pi_1^*\mc E\otimes \Omega^1_{P_1}\
&
\xrightarrow{\pi_1^*\Psi} &
\pi_1^*\mc E\otimes\Omega^2_{P_1}
\end{array}
\end{equation}
It is straightforward to check that
$q \otimes 1 \circ \pi_1^*\Psi$ factorizes through
$\Omega^1_{P_1}(1)$, thereby defining the dashed arrow $a$,
and that for $\kappa\in\mc K^1$
we have $a(\kappa\otimes t_i)=d\kappa\otimes t_i\in(\Omega^2_{P_1}/\mc
K^1\wedge\Omega^1)(1)$. Thus
$M'{}^2=\Omega^2_{P_1}/(\mc K^1\wedge\Omega^1+d\mc K^1)$ is obtained by
factoring out on the upper
right of
\eqref{eqn4.4.5} by the image of the composition across the top twisted by $\mc
O_{P_1}(-1)$. Note
that the composition across the bottom is the curvature $F(\Psi)$. If we write
$(f_{ij})$ for the
curvature matrix with respect to the basis $e_0,e_1,\ldots$, we find using
\eqref{eqn4.4.4} that
e.g. on the open set $t_0\ne 0$, $d\mc K^1$ is generated by elements
\begin{gather}\label{eqn4.4.6}t_0^{-1}(q\otimes 1)F(\Psi)(e_j-(t_j/t_0)e_0) =
\\
\notag \sum_if_{ij}(t_i/t_0)-(t_j/t_0)\sum_kf_{k0}(t_k/t_0)
\end{gather}
On the other hand, the map $\iota$ in \eqref{eqn4.3.2} is given by
\begin{gather}\mc O_{P_1}(-1)\hookrightarrow \pi_1^*\mc E^\vee;\
t_i^{-1}\mapsto\sum_{j}(t_j/t_i)e_j^\vee \notag \\
\label{eqn4.4.7}\Omega^1_{P_1/T}\hookrightarrow\pi_1^*\mc
E(-1)\hookrightarrow\pi_1^*\mc E\otimes\pi_1^*\mc E^\vee \\
d(t_j/t_0)\mapsto (e_j-(t_j/t_0)e_0)\otimes
\sum_i(t_i/t_0)e_i^\vee \notag
\end{gather}
The map $\mu$ from \eqref{eqn4.3.2} is given by $\mu(e_i^\vee\otimes e_j) =
f_{ij}$ hence by
\eqref{eqn4.4.7} we get
\begin{equation}\label{eqn4.4.8}\mu(d(t_j/t_0))=
\sum_i(t_i/t_0)f_{ij}-\sum_i(t_it_j/t_0^2)f_{i0} .
\end{equation}
Comparing \eqref{eqn4.4.6} and \eqref{eqn4.4.8}, we conclude that
\eqref{eqn4.4.3} holds for $i=1$.
The inductive step is precisely the same. We have $P_{i+1} = \mb P(\mc G_{i})$
for some subbundle
$\mc G_{i}\subset\mc G_{i-1}|P_{i}$. The question is local, so we may assume
$\mc G_i$ is free. We
assume inductively that $G_{i-1}$ has a $M_{i-1}=\Omega^*_{P_{i-1}}/\mc
I_{i-1}^*$-connection. Define
$\tilde {M_{i-1}}=\Omega^*_{P_i}/\mc I^*_{i-1}\cdot\Omega^*_{P_i}$, so $\mc
G_{i-1}|P_i$ has a
$\tilde
M_{i-1}$-connection. One factors out by the image $\mc K^1_i$ of
$\Omega^1_{P_i/P_{i-1}}$ as in
\eqref{eqn4.4.2} to define $M^1_i$ and the writes down a diagram like
\eqref{eqn4.4.5} to compare
$d\mc K^1_i$ with the image of $\mu$ as in \eqref{eqn4.3.2}.
At this point it is good to  remark that the curvature
$F_{\tau}: \sO_Q(1) \to \pi^*\Omega^1_T \otimes \sO_Q(1) \to M^2
\otimes \sO_Q(1)$ does not vanish. For example, for $N=2$, one
has $F_{\tau}(t_0) = (f_{00} + f_{01}(t_1/t_0))t_0$.

The remaining assertions in proposition
\ref{prop4.4.1} are easily verified.
\end{proof}

\refstepcounter{subsection}
\begin{prop}\label{prop4.5.1} We have $\mathbb R\pi_*M^i \cong \pi_*M^i$ for
$i< q$. The complex
$H^0(T,\pi_*M^*)$ has no cohomology in odd degrees $< q$. For $2n<q$, the map
$$\delta : \Bbb H^{n-1}(P,M^n\to \ldots\to M^{2n-1}) \to H^0(P,M^{2n})
$$
is injective.
\end{prop}

\subsection{} We postpone the proof of proposition \ref{prop4.5.1} for a while
in order to finish
the proof of theorem \ref{thm4.0.1}. Note first that since the curvature of the
original bundle $E$
on
$X$ is zero, the construction of proposition \ref{prop4.4.1}
above applied to $E$ and the flag bundle
$Q$ yields a structure of differential graded algebra on
$\pi^*\Omega^*_X$, and we have (from
\eqref{eqn4.1.1}) a pullback map of complexes of sheaves on $P$
\begin{equation}\varphi^* : M^* \to R\varphi_*\pi^*\Omega^*_X
\end{equation}
coming from $\varphi^* M^i \to \pi^*\Omega^i_X$.

We will construct classes $\tilde c$ and $\tilde w$ in $\mathbb
H^{n-1}(P,M^n\to \ldots\to
M^{2n-1})$ such that with reference to the maps
\begin{equation}\label{eqn4.6.2}\minCDarrowwidth2pt \begin{CD}
\mathbb H^{n-1}(Q,\pi^*\Omega^n_X\to\ldots\to\pi^*\Omega^{2n-1}_X)\<
\gamma <
\cong <
H^0(X,\Omega^{2n-1}_X)/dH^0(X,\Omega^{2n-2}_X)\\
\V \alpha V V \\
\mathbb H^{n}(Q,\mc {K}^m_n\to\pi^*\Omega^n_X\to\ldots\to\pi^*\Omega^{2n-1}_X)
\\
\A \beta A A  \\
\mathbb H^{n}(Q,\mc {K}^m_n\to \pi^*\Omega^{{\ge n}}_{X})
\end{CD}
\end{equation}
we have
\begin{align}
\beta\pi^*(c_n(E,\nabla)) &= \alpha\varphi^*\tilde c\notag \\
w_n(E,\nabla) &= \gamma^{-1}\varphi^*\tilde w\label{eqn4.6.3}
\end{align}
(Note that to avoid confusion between
$H^0(\Omega^{2n-1}_X/d\Omega^{2n-2}_X)$ and \linebreak
$H^0(\Omega^{2n-1}_X)/dH^0(\Omega^{2n-2}_X)$,
it is a good idea here to localize more and replace $X$
by its function field $\text{Spec}(k(X))$. Note also
that $\beta$ is always injective, and that $\gamma$ is an isomorphism
because $X$ is affine).

We then show
\begin{equation}\label{eqn4.6.4} \delta\tilde c = \delta\tilde w\in
H^0(P,M^{2n}),
\end{equation}
whence, by proposition \ref{prop4.5.1} we have $\tilde c=\tilde w$. Now
consider the analogue of
\eqref{eqn4.6.2} down on $X$, with $\pi^*\Omega$ replaced by $\Omega$. Write
$\alpha_X,\beta_X,\gamma_X$ for the corresponding maps. The assertion of
theorem \ref{thm4.0.1} is
\begin{equation}\label{eqn4.6.5}\gamma_X(w_n(E,\nabla)) =
\beta_X(c_n(E,\nabla)).
\end{equation}
It follows from \eqref{eqn4.6.2} and evident functoriality of $\pi^*$ that
\eqref{eqn4.6.5} holds
after pullback by $\pi^*$. Theorem \ref{thm4.0.1} then follows from

\begin{lem}\label{lem4.6.1} The pullback
\begin{multline*}\pi^* :\mathbb H^{n}(X,\mc
K^m_n\to\Omega^n_X\to\ldots\to\Omega^{2n-1}_X)\to \\
\mathbb H^{n}(Q,\mc K^m_n\to\pi^*\Omega^n_X\to\ldots\to\pi^*\Omega^{2n-1}_X)
\end{multline*}
is injective.
\end{lem}
\begin{proof}[Proof of lemma.] This is central to the splitting principle
involved in the
construction of characteristic classes in \cite{EII}. (See the argument on p.
52 in the proof of
theorem 1.7 of (op. cit.).) An evident diagram chase involving cohomology of
the $\mc K$-sheaves and
the sheaves $\Omega$ and $\pi^*\Omega$ reduces one to showing the $\mc
K$-cohomology groups
$H^{n-1}(X,\mc K^m_n)$ and $H^{n-1}(Q,\mc K^m_n)$ have the same image in
$$\H^{n-1}(Q,\pi^*\Omega^n_X\to\ldots\to\pi^*\Omega^{2n-1}_X).$$
This follows because the
multiplication map
$$H^{n-2}(Q,\mc K^m_{n-1})\otimes H^1(Q,\mc K_1) \to H^{n-1}(Q,\mc
K^m_n)/\pi^*H^{n-1}(X,\mc K^m_n)
$$
is surjective. The line bundles on $Q$ have integrable $\tau$-connections in
the sense of
\cite{EII}, so their classes in $\H^1(Q,\pi^*\Omega^1_X\to\ldots)$ vanish.
\end{proof}

\subsection{} We turn now to the construction of the classes $\tilde c$ and
$\tilde w$. One has
\begin{multline*} w(\mc E,\Psi,P_n)\in
H^0(\Omega^{2n-1}_T)/dH^0(\Omega^{2n-2}_T)\\
 \cong \mathbb H^{n-1}(T,\Omega^n_T\to \ldots\to \Omega^{2n-1}_T)
\end{multline*}
We define $\tilde w$ by the natural pullback
\begin{equation} \tilde w = \pi^*w(\mc E,\Psi,P_n)\in \mathbb H^{n-1}(P,M^n\to
\ldots\to
M^{2n-1}).
\end{equation}
It follows that
\begin{equation}\label{eqn4.7.3} \delta(\tilde w) = \pi^*(dw(\mc E,\Psi,P_n)) =
\pi^*(P_n(F(\Psi)))\in H^0(P,M^{2n}).
\end{equation}

To construct $\tilde c$, we remark first that the map
$$\mb H^{n-1}(P,M^n\to\ldots\to M^{2n-1})\to \mb H^{n}(P,\mc K^m_n\to
M^n\to\ldots\to M^{2n-1})
$$
is injective, so it suffices to construct
\begin{equation}\label{eqn4.7.4}\tilde c\in\ker\bigl(\mb H^{n}(P,\mc K^m_n\to
M^n\to\ldots\to
M^{2n-1})\to H^{n}(P,\mc K^m_n)\bigr).
\end{equation}
This injectivity follows easily from the structure of $H^{n-1}(P,\mc K^m_n)$,
given that $P$ is a
flag bundle over affine space, but in fact the construction of $\tilde c$ as in
\eqref{eqn4.7.4}
would suffice for our purposes anyway, so we won't give the argument in detail.

Let $\ell_i$ be the rank one subquotients of
$\pi^*\mc E$. A basic result from \cite{EI} is that $\pi^*\mc E$ admits a
``connection" with
values in $M^*$,
\begin{equation}\pi^*\mc E \to \pi^*\mc E\otimes_{\mc O_P}M^1
\end{equation}
and that the
filtration defining the
$\ell_i$ is horizontal for this ``connection". Thus there exist local
transition functions
$f_{\alpha,\beta}^i$ and local connection forms $\omega^i_{\alpha}\in M^1$
verifying
\begin{equation}
d {\rm log} f^i = \partial \omega^i ,
\end{equation}
and thus defining $\ell_i \in {\Bbb H}^1(P, {\sK}_1 \to
M^1)$. Here $\partial$ is the Cech differential. Then
$\tilde c$ is defined by the cocyle
\begin{equation}
(x', x^n, ...,x^{2n-1}) \in ({\mc C}^n({\sK}_n) \times
{\sC}^{n-1}(M^n) ...\times {\sC}^0(M^{2n-1}))_{d - \partial}
\end{equation}
with
\begin{gather}
x'= \sum_{i_1<....<i_n} f^{i_1} \cup....\cup f^{i_n} \notag \\
x^n = \sum_{i_1<....<i_n}\omega^{i_1}\wedge \partial \omega^{i_2} \wedge
...\wedge \partial \omega^{i_n} \notag \\
x^{n+1} = \sum_{i_1<....<i_n} \delta\omega^{i_1} \wedge\omega^{i_2}
\wedge \partial\omega^{i_3} \wedge ....\wedge \partial \omega^{i_n} \\
....\notag \\
x^{2n-1} = \sum_{i_1<....<i_n} \delta \omega^{i_1}\wedge . ...\wedge\delta
\omega^{i_{n-1}} \wedge \omega^{i_n} .\notag
\end{gather}
The cup products ``$\cup$" here are Cech products.
By definition \cite{EII},
$\beta\pi^*(c_n(E, \nabla)=  \varphi^*\tilde{c}$.
Applying $\delta$ to the last equation, it follows that the image of $\tilde c$
in $H^0(P,
M^{2n})$ is
\begin{equation}
\sum_{i_1<....<i_n} F(\ell_{i_1}) \wedge ....\wedge F(\ell_{i_n})
{}.
\end{equation}
This is exactly $P_n(F(\oplus\ell_i)) = \pi^*P_n(F(\Psi))$.
(As $M^*$ is a quotient complex of $\Omega^*_P$ by
\ref{prop4.4.1}, (iii),
invariance for $P_n$ guarantees
independence of the choice of local bases for $\pi^*\mc E$.)
Comparing this with \eqref{eqn4.7.3}, we
conclude $\delta\tilde c=\delta\tilde w$ so \eqref{eqn4.6.4} holds.

\subsection{} We turn now to proof of proposition \ref{prop4.5.1}.
\begin{prop}\label{prop4.8.1} The Koszul complex associated to \eqref{eqn4.3.3}
\begin{gather}\label{eqn4.8.1}\ldots \to
\Omega^2_{P/T}\otimes_{\mathcal O_P}\pi^*\Omega^*_T[-4]\to
\Omega^1_{P/T}\otimes_{\mathcal
O_P}\pi^*\Omega^*_T[-2]\xrightarrow{\mu\otimes 1} \notag \\
 \pi^*\Omega^*_T \to M^* \to 0
\end{gather}
is exact in degrees $< q$.
\end{prop}

To clarify and simplify the argument, we will use commutative algebra. Let $B$
be a commutative ring.
Let $C$ be a commutative, graded $B$-algebra, and let $S$ be a graded
$C$-module. Let $Z$ be a
finitely generated free $B$-module with generators $\epsilon_\alpha$, and let
$\nu : Z \to C_2$, with
$\nu(\epsilon_\alpha)=f_\alpha$. Let $\mc I\subset C$ be the ideal generated by
the $f_\alpha$. Write
$\text{gr}_{\mc I}(S) := \oplus \mc I^nS/\mc I^{n+1}S$. Note $\text{gr}_\mc
I(S)$ is a graded module for the
symmetric algebra $B[Z]$ (with $Z$ in degree $2$). The dictionary we have in
mind is
\begin{align}\label{eqn4.8.2} B =& \Gamma(T,\mc O_T) \notag \\
C =& \Omega^{\text{even}}_T\subset S = \Omega^*_T \notag \\
Z =& Hom(\mc E,\mc E) \\
f_\alpha = f_{ij} =&\ \text{entries of curvature matrix} \notag
\end{align}

\begin{lem}\label{lem4.8.2} Let $d\ge 2$ be given. The following are
equivalent.
\begin{enumerate}
\item[(i)] The evident map
$$\rho : (S/\mc IS)[Z] := (S/\mc IS)\otimes_BB[Z] \to \text{gr}_\mc I(S)
$$
is an isomorphism in degrees $\le d$.
\item[(ii)] For all $\alpha$, the multiplication map
\begin{equation}\label{eqn4.8.3} f_\alpha : S/(f_1,\ldots,f_{\alpha-1})S \to
S/(f_1,\ldots,f_{\alpha-1})S
\end{equation}
is injective in degrees $\le d$.
\end{enumerate}
\end{lem}
\begin{proof} This amounts to redoing the argument in Chapter 0, \S
(15.1.1)-(15.1.9) of \cite{G} in
a graded situation, where the hypotheses and conclusions are asserted to hold
only in degrees $\le
d$. The argument may be sketched as follows.

\noindent {\bf Step 1.} Suppose $\alpha=1$, and write $f=f_1$. Let
$\text{gr}(S)=\oplus
f^nS/f^{n+1}S$. Suppose the kernel of multiplication by $f$ on $S$ is contained
in degrees $> d$.
Then the natural map $ \varphi : (S/fS)[T] \to \text{gr}(S)$ is an isomorphism
in degrees $\le d$.
Here, of course, $T$ is given degree = degree$(f) = 2$.

Indeed, $\varphi$ is always surjective, and injectivity in degrees $\le d$
amounts to the assertion
that for $x\in S$ of degree $\le d-2k$ with $f^kx=f^{k+1}y$, we have $x=fy$.
This is clear.

\noindent {\bf Step 2.} Suppose now the condition in (ii) holds. We prove (i)
by induction on
$\alpha$. We may assume by step 1 that $\alpha > 1$. Let $\mc J$ (resp. $\mc
I$) be the ideal
generated by
$f_1,\ldots, f_{\alpha-1}$ (resp. $f_1,\ldots, f_\alpha$). Write
$\text{gr}_{\mc J}(S) = \oplus
J^nS/J^{n+1}S$. By induction, we may assume
\begin{equation}\label{eqn4.8.4}S/\mc JS[T_1,\ldots,T_{\alpha-1}]
\to\text{gr}_{\mc J}(S)
\end{equation}
is an isomorphism in degrees $\le d$. We have to show the same for
\begin{equation}\label{eqn4.8.5}\psi:\text{gr}_{\mc
J}(S)/f_\alpha\text{gr}_{\mc J}(S)[T_\alpha] \to
\text{gr}_{\mc I}(S).
\end{equation}
By \eqref{eqn4.8.4} we have that multiplication by $f_\alpha$ on
$\text{gr}_{\mc J}(S)$ is injective
in degrees $\le d$ (where the degree grading comes from $S$, not the
$\text{gr}_\mc J$ grading). An
easy argument shows the multiplication map
\begin{equation}\label{eqn4.8.6} f_\alpha : S/\mc J^{r}S\hookrightarrow S/\mc
J^{r}S
\end{equation}
is injective in degrees $\le d$ for all $r$.
 Define
\begin{gather}\label{eqn4.8.7} (Q_k)_i = \sum_{j\le k-i}(\text{gr}_{\mc
J}^{k-j}(S)/f_\alpha\text{gr}_{\mc J}^{k-j}(S))T^j \notag \\
(Q_k)_0 = Q_k \qquad (Q_k)_{k+1}=(0) \\
\text{gr}^i(Q_k) = (\text{gr}_{\mc
J}^{k-i}(S)/f_\alpha\text{gr}_{\mc J}^{k-i}(S))T^i\notag
\end{gather}
Define
$$Q_k' = \psi(Q_k)\qquad (Q_k')_i = \psi((Q_k)_i)\qquad
\text{gr}^i(Q_k')=(Q_k')_i/(Q_k')_{i+1}.
$$
The map $\psi$ is surjective, so it will suffice to show the maps
\begin{equation}\label{eqn4.8.8}\text{gr}^i(Q_k) \to \text{gr}^i(Q_k')
\end{equation}
are injective in $(S)$-degrees $\le d$. The left hand side is
$$\mc J^iS/\bigl(f_\alpha\mc J^iS+\mc J^{i+1}S\bigr) T^{k-i}.
$$
The right hand side of \eqref{eqn4.8.8} is the image of
$$\mc J^kS+f_\alpha\mc J^{k-1}S+\ldots +
f_\alpha^{k-i-1}\mc J^{i+1}S$$ in $\mc I^kS/\mc I^{k+1}S.$
What we have to show is that for $x\in
\mc J^iS$ of degree $\le d-2(k-i)$, the inclusion
\begin{equation}\label{eqn4.8.9}f_\alpha^{k-i}x\in \mc J^kS+f_\alpha\mc
J^{k-1}+\ldots+f_\alpha^{k-i-1}\mc J^{i+1}+\mc I^{k+1}S
\end{equation}
implies $x\in f_\alpha\mc J^iS+\mc J^{i+1}S$. The right side of
\eqref{eqn4.8.9} is contained in $\mc
J^{i+1}S+\mc I^{k+1}S\subset\mc J^{i+1}S+f_\alpha^{k+1-i}S$. Multiplication by
$f_\alpha$ on $S/\mc
J^{i+1}S$ is injective in degrees $\le d$ by \eqref{eqn4.8.6}, so
$f_\alpha^{k-i}x\in
f_\alpha^{k+1-i}S+\mc J^{i+1}S$ implies there exists $y\in S$ such that
$x-f_\alpha y\in \mc
J^{i+1}S$. Since $x\in \mc J^iS$, we have $f_\alpha y\in \mc J^iS$ whence by
\eqref{eqn4.8.6} again,
$y\in \mc J^iS$ so $x\in f_\alpha\mc J^iS + \mc J^{i+1}S$. This completes the
verification of step 2.

\noindent {\bf Step 3.} It remains to show $(i)\Rightarrow (ii)$. Again we
argue by induction on
$\alpha$. Suppose first $\alpha=1$. Given $x\in S$ non-zero of degree $\le d$
such that $f_1x=0$, it
would follow from (i) that $x\in f_1^NS$ for all $N$, which is
ridiculous by reason of degree. Now
suppose $\alpha \ge 2$ and that (i) implies (ii) for $\alpha-1$. By assumption
the map
\begin{equation}\label{eqn4.8.10} S/\mc IS[T_1,\ldots,T_\alpha] \to
\text{gr}_{\mc I}(S)
\end{equation}
is an isomorphism in degrees $\le d$. In particular, multiplication by $f_1$ is
injective in
degrees $\le d$ on $\text{gr}_{\mc I}S$. Arguing as above, an $x\in S$ of
degree $\le d$ such that
$f_1 x=0$ would lie in $\mc I^NS$ for all $N$, a contradiction. Thus the first
step in (ii)
holds. To finish the argument, we may factor out by $f_1$, writing $\Bar S =
S/f_1 S$.
Let $\mc K$ be the ideal generated by $f_2,\ldots,f_\alpha$. Factoring out by
$T_1$ on both sides of
\eqref{eqn4.8.10} yields
$$\Bar S/\mc K\Bar S[T_2,\ldots,T_{\alpha}] \to \text{gr}_{\mc K}(\Bar S)
$$
injective in degrees $\le d$. We conclude by induction that (ii) holds for
$\Bar S$.
\end{proof}

Continuing the dictionary from \eqref{eqn4.8.2} above, the ring $R$ and the
module $W$ in the lemma
below correspond to the ring of functions on some affine in $P$ and the module
of $1$-forms
$\Omega^1_{P/T}\subset \pi^*Hom(\mc E,\mc E)$.

\begin{lem}\label{lem4.8.3} Let notation be as above, and assume $\nu : Z\to
C_2$ satisfies the
equivalent conditions of lemma \ref{lem4.8.2}. Let $R$ be a flat $B$-algebra,
and let $W\subset
Z\otimes_BR$ be a free, split $R$-submodule with basis $g_\beta$. Then the
multiplication maps
$$g_\beta : S\otimes_BR/(g_1,\ldots,g_{\beta-1})S\otimes_BR \to
S\otimes_BR/(g_1,\ldots,g_{\beta-1})S\otimes_BR
$$
are injective in degrees $\le d$.
\end{lem}
\begin{proof} Assume not. We can localize at some prime of $R$ contained in the
support for some
element in the kernel of multiplication by $g_\beta$ and reduce to the case $R$
local. Then we may
extend $\{g_\beta\}$ to a basis of $Z\otimes_BR$ and use the implication
(i)$\implies $(ii) from
lemma \ref{lem4.8.2}. Note that $\nu\otimes 1 : Z\otimes R \to C_2\otimes R$
satisfies (i) by
flatness.
\end{proof}

\bigskip

\begin{lem}\label{lem4.8.4} With notations as above, assume $Z$ satisfies the
conditions of lemma
\ref{lem4.8.2} for some $d \ge 2$. Let $J\subset R\otimes_BC$ be the ideal
generated by $(1\otimes
\nu)(W)$. Then the Koszul complex
\begin{multline*}\ldots \to \wedge^2W\otimes_R(R\otimes_BS) \to
W\otimes_R(R\otimes_B S) \to \\
R\otimes_BS \to (R\otimes_BS)/J
\to 0
\end{multline*}
is exact in degrees $\le d$.
\end{lem}
\begin{proof}To simplify notation, let $A=R\otimes_BC,\ M=R\otimes_BS,\
V=W\otimes_RC$, so the
Koszul complex becomes
$$\ldots\wedge^2V\otimes_AM\to V\otimes_AM\to M \to M/JM \to 0.
$$
We argue by induction on the rank of $V$. If this rank is $1$, the assertion is
that the sequence
$$0 \to M \xrightarrow{g_1} M \to M/g_1M \to 0
$$
is exact in degrees $\le d$, which follows from lemma \ref{lem4.8.3}. In
general, if $V$ has an
$A$-basis $g_1,\ldots,g_\beta$, let $V'$ be the span of
$g_1,\ldots,g_{\beta-1}$. By induction, the
Koszul complex
$$\ldots\wedge^2V'\otimes M \to V'\otimes M \to M
$$
is a resolution of $M/(g_1,\ldots,g_{\beta-1})M$ in degrees $\le d$. If we
tensor this module with
the two-term complex $A\xrightarrow{g_\beta}A$ we obtain a complex which by
lemma \ref{lem4.8.3} is
quasi-isomorphic to $M/(g_1,\ldots,g_\beta)M$ in degrees $\le d$. On the other
hand, this complex is
quasi-isomorphic to the complex obtained by tensoring $A\xrightarrow{g_\beta}A$
with the
above $V'$-Koszul complex, and this tensor product is identified with the
$V$-Koszul complex.
\end{proof}

For our application, $B=\mathbb C[x_{ij}^{(k)},y_{ij}^{(k)}]$ is the polynomial
ring in two sets of
variables, with $1\le i,j\le N=\dim(E)$ and $1\le k\le q$ for some large
integer $q$. Let $\Omega$
be the free $B$-module on symbols $dx_{ij}^{(k)}$ and $dy_{ij}^{(k)}$. Let
$S=\bigwedge_B\Omega$,
graded in the obvious way with $dx$ and $dy$ having degree 1, and let
$C=S_{\text{even}}$ be the
elements of even degree. Define
\begin{equation} a_{ij} = \sum_{\ell=1}^q x^{(\ell)}_{ij}dy^{(\ell)}_{ij};\quad
f_{ij}=da_{ij}-\sum_{m=1}^Na_{im}\wedge a_{mj}.
\end{equation}
We have
\begin{equation}f_{ij} = \sum_{\ell=1}^q dx^{(\ell)}_{ij}dy^{(\ell)}_{ij}-
\sum_{m,\ell,p}x^{(\ell)}_{im}x^{(p)}_{mj}dy^{(\ell)}_{im}\wedge dy^{(p)}_{mj}.
\end{equation}
Now give $S$ and $C$ a second grading according to the number of $dx$'s in a
monomial. We denote
this grading by $z = \sum z(j)$. For example, $f_{ij} = f_{ij}(1)+f_{ij}(0)$
with
$f_{ij}(1)=\sum_kdx_{ij}^{(k)}\wedge dy_{ij}^{(k)}$. Let $Z$ be the free
$B$-module on symbols
$\epsilon_{ij}$, with $1\le i,j\le N$. We consider maps
$$\mu, \mu(1) : Z \to C_2;\ \mu(\epsilon_{ij}) = f_{ij};\
\mu(1)(\epsilon_{ij})=f_{ij}(1)
$$

\begin{lem} Suppose the map $\mu(1)$ above satisfies the conditions of lemma
\ref{lem4.8.2} above for
some $d\ge 2$. Then so does $\mu$.
\end{lem}
\begin{proof} Suppose
$$f_\alpha\ell_\alpha = \sum_{1\le \beta \le \alpha-1}f_\beta\ell_\beta
$$
with the $\ell_\beta$ homogeneous of some degree $< d$. Write
$$\ell_\beta = \sum_{0\le j\le r}\ell_\beta(j);\quad 1\le \beta \le \alpha
$$
such that $\ell_\beta(r)\ne 0$ for some $\beta$. We have
\begin{equation}\label{eqn4.8.14} f_\alpha(1)\ell_\alpha(r) = \sum_{1\le
\beta\le \alpha-1}
f_\beta(1)\ell_\beta(r)
\end{equation}
We want to show $\ell_\alpha$ belongs to the submodule generated by
$f_1,\ldots,f_{\alpha-1}$, and
we will argue by double induction on $r$ and on the set
$$\mathcal A = \{\beta\le \alpha\ |\
\ell_\beta(r)\ne 0\}.
$$
If
$r=0$ and
$\ell_\alpha\ne 0$, we get a contradiction from
\eqref{eqn5}, since we have assumed the $\ell_\beta$ have degree $< d$, and
$\ell_\alpha(0)$
cannot lie in the ideal generated by the $f_\beta(1)$. Assume now $r\ge 1$.

\noindent Case 1.- Suppose $\ell_\alpha(r)\ne 0$. From the above, we conclude
we can write
$$\ell_\alpha(r) = \sum_{\beta\in\mathcal A, \beta\ne
\alpha}m_\beta(r-1)f_\beta(1)
$$
Define
\begin{align}\label{eqn4.8.15} \ell_\alpha'= & \ell_\alpha
-\sum_{\beta\in\mathcal A, \beta\ne
\alpha}m_\beta(r-1)f_\beta \\
\ell_\beta' = & \ell_\beta - m_\beta(r-1)f_\alpha\ ;\quad \beta\in\mathcal A,\
\beta\ne \alpha.
\notag
\end{align}
We have still, taking $\ell_\beta' = \ell_\beta$ for $\beta\notin \mathcal A$
\begin{equation}\label{eqn4.8.16} f_\alpha\ell_\alpha'= \sum_{1\le \beta\le
\alpha-1}f_\beta\ell_\beta'
\end{equation}
Since $\ell_\beta'(r)=\ell_\beta(r)=0$ for $\beta\notin \mathcal A$,
$\ell_\alpha'(r)=0$, and
$\ell_\beta'(s)=0$ for $s>r$ and all $\beta$; the inductive hypothesis says
$\ell_\alpha'$ lies in
the ideal generated by the
$f_\beta$ for
$\beta<\alpha$. It follows from \eqref{eqn4.3.2} that $\ell_\alpha$ lies in
this ideal also.

\noindent Case 2.- $\ell_\alpha(r)=0$. Choose $\gamma\in\mathcal A$. We have
$$\sum_{\beta\in\mathcal A}f_\beta(1)\ell_\beta(r)=0.
$$
Since the $f_\beta(1)$ are assumed to satisfy the equivalent hypotheses of
lemma \ref{lem4.6.1}, we
can write
$$\ell_\gamma(r) = \sum_{\beta\in\mathcal A,\beta\ne
\gamma}m_\beta(r-1)f_\beta(1).$$
As in \eqref{eqn4.3.2}, we write
\begin{align}\label{eqn4.8.17} \ell_\gamma'= & \ell_\gamma
-\sum_{\beta\in\mathcal A, \beta\ne
\gamma}m_\beta(r-1)f_\beta \\
\ell_\beta' = & \ell_\beta + m_\beta(r-1)f_\gamma\ ;\quad \beta\in\mathcal A,\
\beta\ne \gamma.
\notag
\end{align}
Again, taking $\ell_\beta'=\ell_\beta$ for $\beta\notin \mathcal A$, we get
\eqref{eqn4.3.3}, so we
may conclude by induction.
\end{proof}

\begin{lem}The map $\mu(1)$ defined by
\begin{equation}\label{eqn4.8.18}
\mu(1)(\epsilon_{ij})=f_{ij}(1)=\sum_{\ell=1}^q
dx^{(\ell)}_{ij}\wedge dy^{(\ell)}_{ij}
\end{equation}
satisfies the hypotheses of lemma \ref{lem4.6.1} with $d=q-1$.
\end{lem}
\begin{proof}Let $V_{ij}$ be the $\mathbb C$-vector space of dimension $2q$
with basis the
$dx^{(\ell)}_{ij}$ and the $dy^{(\ell)}_{ij}$. Write $V=\oplus V_{ij}$. We have
$$ S = \bigwedge V\otimes B \cong \otimes_{i,j}(\bigwedge V_{ij})\otimes B.
$$
It is convenient to well-order the pairs $ij$, writing
$f_\alpha(1)=f_{ij}(1)\in \bigwedge
V_\alpha$. We have
\begin{multline*} S/(f_1,\ldots,f_{\alpha-1})S \cong \\
\bigotimes_{\beta\ge \alpha}\bigl(\bigwedge
V_\beta\bigr)\otimes\bigotimes_{\beta <
\alpha}\biggl(\bigl(\bigwedge V_{\beta}\bigr)/\bigl(f_\beta(1)\bigwedge
V_{\beta}\bigr)\biggr)\otimes B.
\end{multline*}
It is clear from this that multiplication by $f_\alpha(1)$ will be injective in
a given degree $d$
if the multiplication map $f_\alpha(1) : \bigwedge V_\alpha \to \bigwedge
V_\alpha$ is injective in
degrees $\le d$. It is clear from the shape of $f_\alpha(1)$ in
\eqref{eqn4.6.2} that multiplication
by
$f_\alpha(1)$ will be injective in degrees $\le q-1$.
\end{proof}

This completes the proof of proposition \ref{prop4.8.1} above.
\refstepcounter{subsection}
\begin{prop}\label{prop4.9.1} We have for $i< q$
\begin{equation}\label{eqn4.9.1}\pi_*M^0 = \mc O_T;\quad \pi_*M^1 =
\Omega^1_T;\quad R^j\pi_*M^i =
(0);\ j\ge 1
\end{equation}
The sheaf $\pi_*M^i$ admits an increasing filtration
$\text{fil}_\ell(\pi_*M^i),\ \ell\ge 0$ which
is stable under $\delta$ and satisfies
\begin{equation}\label{eqn4.9.2}\text{gr}_j(\pi_*M^i)\cong
H^j(P,\Omega^j_{P/T})\otimes
\Omega^{i-2j}_T\cong CH^j(P)\otimes_{\mathbb Z}\Omega^{i-2j}_T
\end{equation}
for $j\ge 0$. Here $CH^j(P)$ is the Chow group of codimension $j$ algebraic
cycles on $P$. The
differential $\text{gr}_j(\pi_*M^i)\to \text{gr}_j(\pi_*M^{i+1})$ is the
identity on the Chow group
tensored with the exterior derivative on $\Omega^*_T$ up to sign.
\end{prop}

\bigskip

Note that the last assertion in \eqref{eqn4.9.1} implies for $i<q$
$$H^*(P,M^i)\cong \begin{cases}H^0(T,\pi_*M^i)& *=0\\
0& *\ge 1.
\end{cases}
$$
It follows from \eqref{eqn4.9.2} that the complex $H^0(T,M^*)$ has no
cohomology in odd degrees $<
q-1$. (Recall that $T$ has no higher de Rham
cohomology). These assertions imply proposition \ref{prop4.5.1}.

\bigskip

\begin{proof}[Proof of proposition \ref{prop4.9.1}.]
The first two assertions in \eqref{eqn4.9.1} are
clear, because $M^0 = \mc O_P$ and $M^1=\pi^*\Omega^1_T$.
We define
$$G_j = {\rm Im } \Omega^j_{P/T} \otimes \pi^* \Omega^{i-2j}_T
\to \Omega^{j-1}_{P/T} \otimes \pi^* \Omega^{i-2j+2}_T,\  G_0 = M^i
$$
coming from the resolution of $M^i$ in \ref{prop4.8.1}. Then $R^a\pi_*G_j=0$
for $a \neq j$. This proves $R^j\pi_*M^i=0$ for $j \geq 1$. One has
a short exact sequence
$$0 \to R^j\pi_* \Omega^j_{P/T} \otimes \Omega^{i-2j}_T \to R^j\pi_* G_j
\to R^{j+1}\pi_*G_{j+1} \to 0$$
with $R^0\pi_*G_0 = \pi_* M^i$. One defines
\begin{multline*}
\text{fil}_j(\pi_*M^i) = \text{inverse image of} \ R^j\pi_*\Omega^j_{P/T}
\otimes \Omega^{i-2j}_T \\
\text{via} \ R^0\pi_*G_0 \to R^j\pi_*G_j.
\end{multline*}
This proves
\eqref{eqn4.9.2}.

In order to understand the map ${\rm gr}_j(\pi_*M^i) \to
{\rm gr}_j(\pi_*M^{i+1})$, we construct a commutative diagram
\begin{equation}\label{eqn4.9.3}\minCDarrowwidth10pt
\begin{CD}
\> \mu \otimes 1>> \Omega^2_{P/T} \otimes \pi^*\Omega^{i-4}_T \>
\mu \otimes 1>>
\Omega^1_{P/T} \otimes \pi^*\Omega^{i-2}_T \> \mu \otimes 1>>
\pi^*\Omega^i_T \> \mu \otimes 1>>
M^i\\
\noarr
\V V \nabla_{\tau} V  \V V \nabla_{\tau} V  \V V \nabla_{\tau} V  \V V
\nabla_{\tau} V \\
\>\mu \otimes 1>> \Omega^2_{P/T} \otimes \pi^* \Omega^{i-3}_T
\>\mu \otimes 1>>
\Omega^1_{P/T} \otimes \pi^* \Omega^{i-1}_T \>\mu \otimes 1>> \pi^*
\Omega^{i+1}_T \>\mu \otimes 1>> M^{i+1}
\end{CD}
\end{equation}
mapping the resolution of $M^i$ to the resolution of $M^{i+1}$ given
by \ref{prop4.8.1}. To this aim recall that one has an exact
sequence of complexes

\begin{equation}\label{eqn4.9.4}
0 \to K^* \to \Omega^*_P \to M^* \to 0
\end{equation}
with
\begin{equation}\label{eqn4.9.5}
K^i = \Omega^i_{P/T} \oplus \Omega^{i-1}_{P/T} \otimes \pi^*
\Omega^1_T ...\oplus \Omega^1_{P/T} \otimes \pi^* \Omega^{i-1}_T
\oplus \mu(\Omega^1_{P/T}) \wedge \pi^*\Omega^{i-2}_T
\end{equation}
Note that the differential $K^{i-j-1} \to K^{i-j}$ acts as
follows
\begin{gather}\label{eqn4.9.6}
\Omega^j_{P/T} \otimes \pi^* \Omega^{i-1-2j}_T
\to  \\
\notag
\Omega^{j+1}_{P/T} \otimes \pi^* \Omega^{i-1-2j}_T
\oplus \Omega^{j}_{P/T} \otimes \pi^* \Omega^{i-2j}_T
\oplus \Omega^{j-1}_{P/T} \otimes \pi^* \Omega^{i-2j +1}_T.
\end{gather}
To see this, write
\begin{equation}\label{eqn4.9.7}
\Omega^j_{P/T} \otimes \pi^* \Omega^{i-1-2j}_T =
\Omega^j_{P/T} \otimes_{{\sO}_T} \pi^{-1} \Omega^{i-1-2j}_T
\end{equation}
and apply the Leibniz rule with
\begin{equation}\label{eqn4.9.8}
d \Omega^1_{P/T} \subset \Omega^2_{P/T} \oplus \Omega^1_{P/T}
\otimes \pi^* \Omega^1_T \oplus \mu(\Omega^1_{P/T}).
\end{equation}
The corresponding map $\Omega^1_{P/T} \to \mu(\Omega^1_{P/T})$
is of course $\mu$. We denote by $\nabla_{\tau}$ the
corresponding map $\Omega^1_{P/T} \to \Omega^1_{P/T} \otimes
\pi^* \Omega^1_T$ and also by $\nabla_{\tau}$ the induced map
$\Omega^j_{P/T} \otimes \pi^* \Omega_T^{i-1-2j} \to
\Omega^j_{P/T} \otimes \pi^* \Omega_T^{i-2j}$. For $\gamma \in
\Omega^j_{P/T} \otimes \pi^* \Omega_T^{i-1-2j}$, write
$d\gamma= \gamma_{j+1} + \nabla_{\tau}(\gamma) + (\mu \otimes
1)(\gamma)$ with $\gamma_{j+1} \in \Omega^{j+1}_{P/T} \otimes
\pi^* \Omega^{i-1-2j}_T$. The integrability condition $d^2(\gamma) = 0$
in $\Omega^*_P$ says that $ (\mu \otimes 1) \nabla_{\tau}
(\gamma) = \nabla_{\tau}(\mu \otimes 1)(\gamma) \in
\Omega^{j-1}_{P/T} \otimes \pi^* \Omega^{i-j+2}_T$, upto sign.

Thus ${\rm gr}_j\pi_*M^i \to {\rm gr}_j\pi_*M^{i+1}$ is the map
$$R^j\pi_*\nabla_{\tau}: R^j\pi_*\Omega^j_{P/T} \otimes \Omega^{i-2j}_T
\to R^j\pi_*\Omega^j_{P/T} \otimes \Omega^{i-2j+1}_T.$$
Now, $\nabla_{\tau}=
d|K^*$, where $d$ is the differential of $\Omega^*_P$.
Let $\ell_i$ be the rank one subquotients of
$\pi^*\sE$, with local algebraic transition functions $f^i_{\alpha, \beta}$.
Then $R^j\pi_*\Omega^j_{P/T} \otimes \Omega^{i-2j}_T$ is
generated over $\sO_T$
by elements
$\varphi = F \wedge \omega$, with $$ F = d {\rm log}f^{i_1}_{\alpha_0,
\alpha_1} \wedge \cdots \wedge d {\rm log}f^{i_j}_{\alpha_{j-1}, \alpha_j}$$
and $\omega \in \Omega^{i-2j}_T$. Thus $d \varphi = (-1)^j F \wedge d\omega$.
This finishs the proof of the proposition.
\end{proof}

\section{Chern-Simons classes and the Griffiths group}\label{sec:griff}
\subsection{} Our objective in this section is to investigate the vanishing of
the class
$w_n(E,\nabla)$ for a flat bundle $E$ on a smooth, projective variety $X$ over
$\mb C$. We will show
that $w_n=0$ if and only if the
$n$-th chern class $c_n(E)$ vanishes in a ``generalized Griffiths group''
$\text{Griff}^{\,n}(X)$.

Let $X$ be a smooth, quasi-projective variety over $\mb C$. For $Z\subset X$ a
closed subvariety and
$A$ an abelian group, we write $H^*_Z(X,A)$ for the singular cohomology with
supports in $Z$ and
values in $A$. We write
$$H^*_{\mc Z^n}(X,A) = \varinjlim_{Z\subset X\text{ cod. }n}H^*_Z(X,A).
$$
Purity implies that for $Z$ irreducible of codimension $n$,
$$H^p_Z(X,A) = \begin{cases} 0 & p< 2n \\
A(-n) & p=2n
\end{cases}
$$
Here $\mb Z(n) = (2\pi i)^n\mb Z$ and $A(n) := A\otimes \mb Z(n)$. As a
consequence
$$H^p_{\mc Z^n}(X,\mb Z(n)) = \begin{cases} 0 & p< 2n \\
\mc Z^n(X) & p=2n
\end{cases}
$$
where $\mc Z^n(X)$ is the group of codimension $n$ algebraic cycles on $X$.

For $m< n$, define the Chow group of codimension $n$ algebraic cycles modulo
codimension $m$
equivalence by
\begin{equation}CH^n_m(X) := \text{Image}\bigl(\mc Z^n(X) = H^{2n}_{\mc
Z^n}(X,\mb Z(n))\to
H^{2n}_{\mc Z^m}(X,\mb Z(n)\bigr).
\end{equation}
Of course, $CH^*_0(X)$ is the group of cycles modulo homological equivalence.
It follows from
\cite{BO} (7.3) that $CH^n_{n-1}(X)$ is the group of codimension $n$ algebraic
cycles modulo {\bf
algebraic} equivalence.
\begin{defn} The generalized Griffiths group $\text{Griff}^{\,n}(X)$ is defined
to be the kernel of
the map $CH^n_1(X) \to CH^n_0(X)$. In words, the generalized Griffiths group
consists of cycles
homologous to $0$ on $X$ modulo those homologous to $0$ on some divisor in $X$.
\end{defn}
\begin{ex} $\text{Griff}^{\,2}(X)$ is the usual Griffiths group of codimension
$2$ cycles homologous
to zero modulo algebraic equivalence.
\end{ex}

\subsection{} With notation as above, let $\mc H^p(A)$ denote the Zariski sheaf
on $X$ associated to
the presheaf
$U \mapsto H^p(U_{\text{an}},A)$, cohomology for the classical (analytic)
topology with coefficients
in
$A$. The principal object of study in
\cite{BO} was a spectral sequence
\begin{equation}\label{eqn2}E_2^{p,q}(A) = H^p(X_{\text{Zar}}, \mc H^q(A))
\Rightarrow
H^{p+q}(X_{\text{an}},A).
\end{equation}
associated to the ``continuous'' map $X_{\text{an}} \to X_{\text{Zar}}$. This
spectral sequence was
shown to coincide from $E_2$ onward with the ``coniveau'' spectral sequence
\begin{equation}\label{eqn3} E_1^{p,q}(A) = \bigoplus_{x\in \mc Z^p-\mc
Z^{p+1}}H^{q-p}(x,A)\Rightarrow H^{p+q}(X_{\text{an}},A).
\end{equation}
As a consequence of a Gersten resolution for the sheaves $\mc H^p(A)$, one had
\begin{gather}\label{eqn4} H^p(X_{\text{Zar}},
\mc H^q(A)) = (0)\text{ for }p>q \\
\notag H^n(X_{\text{Zar}},\mc H^n(\mb Z(n))) \cong CH^n_{n-1}(X).
\end{gather}
The $E_\infty$-filtration $N^*H^*(X_{\text{an}},A)$ is the filtration by
codimension,
$$N^pH^*(X_{\text{an}},A) = \text{Image}\ (H^*_{\mc Z^p}(X_{\text{an}},A)\to
H^*(X_{\text{an}},A)).
$$
\refstepcounter{subsection}
\begin{prop} With notation as above, there is an exact sequence
\begin{equation}\label{eqn5} 0 \to H^{2n-1}(X_{{\rm an}},\mb Z(n))/N^1 \to
E_n^{0,2n-1} \xrightarrow{d_n} {\rm Griff}^{n}(X) \to 0.
\end{equation}
\end{prop}
\begin{proof} It follows from \eqref{eqn4} that we have
\begin{gather}{}\notag H^{2n-1}(X_{\text{an}},\mb Z(n))\twoheadrightarrow
E_\infty^{0,2n-1}=E_{n+1}^{0,2n-1}\subset E_{n}^{0,2n-1}\subset\cdots  \\
\subset E_2^{0,2n-1} = \Gamma(X,\mc H^{2n-1}(\mb Z(n)))
\end{gather}
and
\begin{equation}CH^n_{n-1}(X)=E_2^{n,n}\twoheadrightarrow
E_3^{n,n}\twoheadrightarrow\cdots
\twoheadrightarrow
E_{n+1}^{n,n}=E_\infty^{n,n}\subset H^{2n}(X_{\text{an}},\mb Z(n)).
\end{equation}
In fact, $E_r^{n,n}\cong CH^n_{n+1-r}(X)$. In particular, $E_n^{n,n}\cong
CH^n_1(X)$. To see this,
one can, for example, use the theory of exact couples \cite{Hu} pp. 232 ff. One
gets an exact
triangle
$$
\begin{array}{lcr} D_r &\makebox[2pt]{$\xrightarrow{i_r}$} & D_r \\
{}_{k_r}\!\nwarrow && \swarrow\!{}_{j_r}\ \\
& \ \makebox[2pt]{$E_r$} &
\end{array}
$$
where in the appropriate degree
$$
D_r= \text{Image}(H^{2n}_{\mc Z^n}(X,\mb Z(n)) \to H^{2n}_{\mc
Z^{n-r+1}}(X,\mb Z(n)))\cong CH^n_{n-r+1}(X),
$$
$i_r = 0$ and $j_r$ is an isomorphism.

The spectral sequence \eqref{eqn2} now yields a diagram with exact rows,
proving the proposition.
\begin{equation}\label{eqn8}\minCDarrowwidth10pt\begin{CD}
0 \>>> H^{2n-1}(X,\mb Z(n))/N^1 \>>>
E_n^{0,2n-1}\>d_n >> {\rm Griff}^{n}(X) \>>>
0 \\
\noarr
 \V V = V \V V = V \V V \cap V \noarr \\
0 \>>> H^{2n-1}(X,\mb Z(n))/N^1 \>>> E_n^{0,2n-1} \> d_n >>
E_n^{n,n} \>>> H^{2n}(X,\mb Z(n))
\end{CD}
\end{equation}
\end{proof}
\refstepcounter{subsection}
\begin{prop}\label{prop2} Let $X$ be smooth and quasi-projective over $\mb C$.
Let $(E,\nabla)$ be a
vector bundle with an integrable connection on $X$.
Let $n\ge 2$ be given, and let $d_n$ be as in
\eqref{eqn8}. Let $c_n(E)$ be the $n$-th Chern class in
$\text{Griff}^{\,n}(X)\otimes\mb Q$. Then
\begin{enumerate}\item[(i)]
$w_n(E,\nabla)\in E_n^{0,2n-1}(\mb C)\subset \Gamma(X,\mc H^{2n-1}(\mb C))$.
\item[(ii)] $d_n(w_n) = c_n(E)$.
\end{enumerate}
\end{prop}
\begin{proof}The spectral sequence \eqref{eqn2} in the case $A=\mb C$ coincides
with the ``second
spectral sequence'' of hypercohomology for
$$H^*(X_{\text{an}},\mb C)\cong \mb H^*(X_{\text{Zar}},\Omega^*_{X/\mb C}).
$$
This is convenient for calculating the differentials in \eqref{eqn2}. Namely,
we consider the
complexes for $m\ge n$
\begin{align}\label{eqn9}\tau_{n,n}\Omega^*:=& \mc H^n(\mb C)[-n]\notag \\
 \tau_{m,n}\Omega^*:= &
\bigl(\Omega^m_X/d\Omega^{m-1}_X\to\ldots\to\Omega^{n-1}_X\to\Omega^n_{\text{closed}}\bigr)[-m];\ m<n.
\end{align}
We have maps
$$\tau_{0,n}\Omega^* \to \tau_{1,n}\Omega^*
\to\cdots\to\tau_{n,n}\Omega^*\to\tau_{n,n+1}\Omega^*\to\cdots\to\tau_{n,\infty}\Omega^*,
$$
and
\begin{gather}
\notag E_r^{0,2n-1}=\text{Image}(H^{2n-1}(X,\tau_{2n-r+1,2n-1}\Omega^*)\to \\
\label{eqn10}          H^{2n-1}(X,\tau_{2n-1,2n-1}\Omega^*)=
\Gamma(X,\mc H^{2n-1})).
\end{gather}
There is a diagram of complexes
\begin{equation}\minCDarrowwidth10pt\label{eqn11}\begin{CD}\bigl[\mc
K^m_n\to\Omega^n\to\cdots\to\Omega^{2n-2}\to\Omega^{2n-1}_{\text{closed}}\bigr]\>a>>
\Omega^\infty\mc K^m_n \\
\V b V V \V c V V \\
\qquad\qquad\tau_{n+1,2n-1}\Omega^*[n-1] \>e>>
\tau_{2n-1,\infty}\Omega^*[n-1],
\end{CD}\end{equation}
where $\Omega^\infty\mc K^m_n$ is the complex $\mc
K^m_n\to\Omega^n\to\Omega^{n+1}\to\cdots$. We
have
\begin{gather}c_n(E,\nabla)\in \mb H^n(X,\Omega^\infty \mc K^m_n)\\
c(c_n(E,\nabla))=w_n(E,\nabla)\in
\mb H^{2n-1}(X,\tau_{2n-1,\infty}\Omega^*)\cong H^0(X,\mc H^{2n-1}).\notag
\end{gather}
 The
map $a$ is the inclusion of a subcomplex, and the quotient has no cohomology
sheaves in degrees
$<n+1$, so $a$ is an isomorphism on hypercohomology in degree $n$.
It follows that $w_n(E,\nabla)$ lies
in the image of the map $e$ in \eqref{eqn11}. By \eqref{eqn10},
this image is $E_n^{0,2n-1}$.

To verify $d_n(w_n)=c_n(E)$, write $\bar\Omega^\infty\mc K^m_n$ for the complex
$$\mc K^m_n\to
\Omega^n/d\Omega^{n-1}\to\Omega^{n+1}\to\cdots,$$
and let $\bar c_n(E,\nabla)\in\mb H^n(X,\bar\Omega^\infty\mc K^m_n)$ be the
image of
$c_n(E,\nabla)$. Consider the distinguished triangle of complexes
\begin{equation}\label{eqn13}\tau_{n,2n-2}\Omega^*[n-1]\to\bar\Omega^\infty\mc
K^m_n
\xrightarrow{\alpha}\mc K^m_n\oplus\tau_{2n-1,\infty}\Omega^*[n-1]
\end{equation}
We have by definition $\alpha(\bar c_n(E,\nabla))=(c_n(E),w_n(E))$, so, writing
$\partial$ for the
boundary map,
\begin{equation}\label{eqn14}\partial(c_n(E))=-\partial(w_n(E,\nabla))\in \mb
H^{2n}(X,\tau_{n,2n-2}\Omega^*).
\end{equation}
Note that the boundary map on $\mc K^m_n$ factors through the dlog map $\mc
K^m_n \to \mc H^n$. Thus
$\partial c_n(E)$ is the image of the Chern class. On the other hand, by
\eqref{eqn4} we have
$$ \mb H^{2n}(X,\Omega^*)\surj \mb H^{2n}(X,\tau_{n,\infty}\Omega^*),
$$
from which it follows by standard spectral sequence theory that the image of
the map
$$H^{2n}(X,\tau_{n,n}\Omega^*)\to\mb H^{2n}(X,\tau_{n,2n-2}\Omega^*)
$$
coincides with $E_n^{n,n}$, and that the boundary map
$$\delta: \Gamma(X,\mc H^{2n-1})\cong\mb H^{2n-1}(X,\tau_{2n-1,\infty}\Omega^*)
\to \mb
H^{2n}(X,\tau_{n,2n-2}\Omega^*)
$$
coincides with $d_n$ from the statement of the proposition on
$\delta^{-1}(E_n^{n,n})=E_n^{0,2n-1}$.
This completes the proof of the proposition.
\end{proof}

Our next objective is to realize the sequence \eqref{eqn5} as an exact sequence
of mixed Hodge
structures. To avoid complications, we replace $\mb Z$ with $\mb Q$ throughout.
More precisely, we
work with filtering direct limits of finite dimensional $\mb Q$-mixed Hodge
structures, where the
transition maps are maps of mixed Hodge structures.
\begin{lem} The spectral sequence \eqref{eqn2} with $A=\mb Q(n)$ can be
interpreted as a spectral
sequence in the category of mixed Hodge structures.
\end{lem}
\begin{proof} The spectral sequence \eqref{eqn3} can be deduced from an exact
couple
(\cite{BO},p.188)
$$\cdots \to H^{p+q}_{\mc Z^p}(X,\mb Q(n))\to H^{p+q}_{\mc Z^{p-1}}(X,\mb
Q(n))\to H^{p+q}_{\mc
Z^{p-1}/\mc Z^p}(X,\mb Q(n))\to\cdots
$$
These groups clearly have infinite dimensional mixed Hodge structures and the
maps are morphisms of
mixed Hodge structures. The lemma follows easily, since \eqref{eqn2} coincides
with the above from
$E_2$ onward.
\end{proof}
\begin{remark} The groups $E_r^{n,n}$ are all quotients of
$$H^{2n}_{\mc Z^n/\mc Z^{n+1}}(X,\mb Q(n))\cong \bigoplus_{z\in X^n}\mb Q
$$
so these groups all have trivial Hodge structures.
\end{remark}
\refstepcounter{subsection}
\begin{prop}\label{prop3} The Chern-Simons class $w_n(E,\nabla)\in
E_n^{0,2n-1}(\mb C)$ lies in
$F^0$ (zeroeth piece of the Hodge filtration) for the Hodge structure defined
by $E_n^{0,2n-1}(\mb
Q(n))$.
\end{prop}
\begin{proof} We have $E_n^{0,2n-1}\subset E_2^{0,2n-1}\subset H^{2n-1}(\mb
C(X),\mb C)$, where the
group on the right is defined as the limit over Zariski open sets. Thus, it
suffices to work ``at
the generic point''. Let $\mc S$ denote the category of triples $(U,Y,\pi)$
with $Y$ smooth and
projective, $\pi : Y\to X$ a birational morphism of schemes, and $U\subset Y$
Zariski open such that
$Y_U$ is a divisor with normal crossings and $U\to \pi(U)$ is an isomorphism.
Using resolution of
singularities, one sees easily that
$$H^n(\text{Spec}(\mb C(X)),\mb C)\cong
\varinjlim_{\mc S}\mb H^n(Y,\Omega^*_{Y}(\log (Y-U))).
$$
The Hodge filtration on the left is induced in the usual way from the first
spectral sequence of
hypercohomology on the right.
\begin{lem}For $\alpha\in\mc S$ let $j_\alpha : U_\alpha \inj Y_\alpha$ be the
inclusion. Then
$$\varinjlim_{\mc S} H^n(Y_\alpha,j_{\alpha*}\mc K^m_{n,U_\alpha}) = (0),\ n\ge
1.
$$
\end{lem}
\begin{proof}[Proof of lemma.] Given $j : U\inj Y$ in $\mc S$ and $z\in
H^n(Y,j_*\mc K^m_{n,U})$,
let $k:\text{Spec}(\mb C(X))\to Y$ be the generic point. We have $H^n(Y,k_*\mc
K^m_{n,\mb C(X)}) =
(0)$ since the sheaf is constant, so there exists $V\subset U$ open of finite
type such that writing
$\ell : V\to Y$, $z$ dies in $H^n(Y,\ell_*\mc K^m_{n,V})$. Let $m: V\to Z$
represent an object of
$\mc S$ with $Z$ dominating $Y$. We have a triangle
\begin{equation}\begin{array}{rcl} H^n(Y,j_*\mc K^m_{n,U}) & \to & H^n(Z,m_*\mc
K^m_{n,V}) \\
\searrow && \nearrow \\
& H^n(Y,\ell_*\mc K^m_{n,V}) &
\end{array}\end{equation}
from which it follows that $z\mapsto 0$ in $\varinjlim_{\mc
S}H^n(Y_\alpha,j_{\alpha,*}\mc
K^m_{n,U_\alpha})$.
\end{proof}

Returning to the proof of proposition \ref{prop3}, write $D_\alpha =
Y_\alpha-U_\alpha$ for
$\alpha\in\mc S$. We see from the lemma that the map labeled
$a$ below is surjective:
\begin{equation}\begin{array}{rcl} c_n(E,\nabla)& \in &
\H^n(X,\Omega^\infty\mc K^m_n) \\
&& \qquad\qquad \downarrow \\
&& \makebox[2in][r]{$\varinjlim_{\mc S}\H^n(Y_\alpha, j_{\alpha,*}\mc
K^m_{n,U_\alpha}\to\Omega^n(\log
(D_\alpha))\to\cdots)$} \\
\nearrow \negmedspace a& &\qquad \qquad \downarrow \\
\makebox[2.1in][l]{$\varinjlim_{\mc S}\H^{n-1}(Y_\alpha, \Omega^n(\log
(D_\alpha))\to\cdots)$}& &\xrightarrow{b} H^{2n-1}(\text{Spec}(\mb C(X)),\mb C)
\end{array}\end{equation}
Since the image of $b$ is $F^0$ for the Hodge filtration on
\linebreak $H^{2n-1}(\text{Spec}(\mb
C(X)),\C)$, and since the composition of vertical arrows maps $c_n(E,\nabla)$
to the
restriction of $w_n(E,\nabla)$ at the generic point, the proposition is proved.
\end{proof}
\refstepcounter{subsection}
\begin{prop}\label{prop4}
The Chern-Simons class $w_n(E,\nabla)\in E_n^{0,2n-1}(\mb Q(n))$.
\end{prop}
\begin{proof}The following diagram is commutative
\begin{equation}\begin{array}{rcl}
\mb H^n(X,\Omega^\infty\mc K^m_n) & \to H^0(X,\mc
H^{2n-1}(\mb C)) \to & H^{2n-1}(\text{Spec}(\mb C(X)),\mb C) \\
\downarrow \varphi(\ref{lem:maps}(3)) && \qquad\downarrow \\
H^{2n-1}(X_{\text{an}},\mb C/\mb
Z(n))&\makebox[1in][r]{$\xrightarrow{\makebox[1.2in]{}}$}&
\makebox[1in][c]{$\ \ \ H^{2n-1}(\text{Spec}(\mb C(X)),\mb C/\mb Z(n))$}.
\end{array}\end{equation}
We know from \cite{EI} and proposition \ref{prop:3.10.1} above that
$$\varphi(c_n(E,\nabla))=c_n^{\text{an}}(E,\nabla) , $$
and, using the deep theorem of Reznikov
\cite{Re}, that this class is torsion. In particular, the image of
$c_n(E,\nabla)$ on the upper
right lies in \linebreak $H^{2n-1}(\text{Spec}(\mb C(X)),\mb Q(n))$.
As a matter of fact, in \ref{thm:5.6.2}, we will only use that
$w_n(E, \nabla) \in E^{0,2n-1}_n(\R(n))$. For this we don't need
the full strength of \cite{Re}, but only that $c^{{\rm an}}_n(E, \nabla)=0
\in H^{2n-1}(X_{{\text an}}, \C/\R(n))$, which is a consequence of
Simpson's theorem \cite{Sim}
asserting that $(E, \nabla)$ deforms to a
$\C$ variation of Hodge structure.
\end{proof}

\begin{thm}\label{thm:5.6.2} Let $X$ be a smooth, projective variety over $\mb
C$. Let $E$ be a vector bundle on $X$,
and let $\nabla$ be an integrable connection on $E$. Then $w_n(E,\nabla)\in
H^0(X,\mc H^{2n-1}(\mb
C))$ vanishes if and only if the cycle class $c_n(E)$ is trivial in
$\text{Griff}^n(X)\otimes\mb Q$.
\end{thm}
\begin{proof} Consider the exact sequence of mixed Hodge structures
\begin{equation} 0 \to H^{2n-1}(X_{\text{an}},\mb Q(n))/N^1 \to
E_n^{0,2n-1}(\mb Q(n)) \to
\text{Griff}^n(X)\otimes\mb Q \to 0.
\end{equation}
Write $H$ for the group on the left. It is pure of weight $-1$, so \linebreak
$H(\mb Q)\cap
F^0H(\mb C) = (0)$. It follows that $w_n(E,\nabla)=0$ if and only if its image
$c_n(E)$ in
$\text{Griff}^n(X)\otimes\mb Q$ vanishes.
\end{proof}

The following corollary is a simple application of the theorem to the example
\eqref{eqn:rk2ex}
discussed in the introduction.
\begin{cor} Let $E,\ X,\ \nabla$ be as above. Assume $E$ has rank $2$, and that
the determinant
bundle is trivial, with the trivial connection. Let $U\subset X$ be affine open
such that $E|U$ is
trivial, and let $\bigl(\begin{smallmatrix}\alpha & \beta \\ \gamma & -\alpha
\end{smallmatrix}\bigr)$ be the connection matrix. Then $c_2(E)
\otimes \Q$ is algebraically equivalent to $0$
on $X$ if and only if there exists a
meromorphic $2$-form $\eta$ on $X$ satisfying $d\eta =
\alpha\wedge d\alpha = \alpha \wedge\beta \wedge \gamma$.
\end{cor}

\section{Logarithmic Poles}\label{sec:log}

In this section we consider a normal crossing
divisor $D \subset X$ on a smooth variety $X$,
the inclusion $j: X-D \to X$,
and a bundle $E$, together with a flat connection
$\nabla: E \to \Omega^1_X({\rm log} D) \otimes E$
with logarithmic poles along $D$. The characteristic
of the ground field $k$ is still 0. Finally recall from
\cite{BO} that one has an exact sequence
\begin{equation}\label{eqn:BO}
0 \to H^0(X, \sH^j) \to H^0(X-D, \sH^j) \rightarrow{\oplus_i{\rm res}}
H^0(k(D_i), \sH^{j-1})
\end{equation}

\refstepcounter{subsection}
\begin{thm}\label{thm:logpoles}Let $(E,\nabla, D)$ be a flat connection with
logarithmic poles. Then
\begin{equation}{}w_n(E, \nabla) \in H^0(X, \sH^{2n-1})
\subset H^0(X-D, \sH^{2n-1}) = H^0(X, j_* \sH^{2n-1}) .
\end{equation}
\end{thm}
\begin{proof}
By \ref{eqn:BO}, one just has to compute the residues of $w_n(E, \nabla)$
along generic points of $D$. So one may assume that the local equation
of $\nabla$ is $A= B \frac{dx}{x} + C $, where $B$ is a matrix
of regular functions, $x$ is the local equation of a smooth component
of $D$, and $C$ is a matrix of regular one forms. Furthermore,
as $dA=A^2=\frac{1}{2}[A,A]$, the formulae of theorem 2.2.1 say that the local
shape of $w_n(E, \nabla)$ is $ {\rm Tr}\lambda A(dA)^{n-1}= \lambda
{\rm Tr} (B \frac{dx}{x} + C)  (dB \frac{dx}{x} + dC)^{n-1}$
for some $\lambda \in \Q$. So upto
coefficient one has to compute
\begin{multline}\label{eqn:basic}{\rm Tr \ Res}  (B \frac{dx}{x} +
C)((dC)^{n-1} +
\sum_{a+b=n-2} (dC)^adB\frac{dx}{x}(dC)^b)= \\
{\rm Tr \ Res}[C \sum_{a+b=n-2} (dC)^adB(dC)^b +
B(dC)^{n-1}]\frac{dx}{x}.
\end{multline}
On the other hand, the integrability condition reads
$$(dB - (C B - B C))
\frac{dx}{x} + d C - C^2 =0,$$
from which one deduces
\begin{gather}dC\frac{dx}{x}=C^2\frac{dx}{x}\label{eqn:1deduct}\\
{\rm   Res} (dB - (CB-BC))\frac{dx}{x} = 0\label{eqn:2deduct}
\end{gather}
Applying \eqref{eqn:1deduct} to \eqref{eqn:basic}, we reduce to calculating
\begin{equation}\label{eqn:2basic}
{\rm Tr \ Res} [\sum_{a+b=n-2} (dC)^a C dB (dC)^b +
B(dC)^{n-1}]\frac{dx}{x}
\end{equation}
 Since we are only interested in calculating \eqref{eqn:2basic} modulo exact
 forms, we can use $d(CB)=dCB-CdB$ and move copies of $dC$ to the right in
\eqref{eqn:2basic} under the trace. The problem becomes to show
\begin{equation}\label{eqn:3basic}{\rm Tr \ Res} B(dC)^{n-1}\frac{dx}{x}
\end{equation}
is exact. It follows from \eqref{eqn:2deduct} that
\begin{equation}
{\rm Tr \ Res }C^{2n-3}dB\frac{dx}{x}={\rm Tr \ Res
}[C^{2n-2}B-C^{2n-3}BC]\frac{dx}{x}
\end{equation}
Bringing the $C$ to the left in the last term changes the sign, so we get by
\eqref{eqn:1deduct}
\begin{equation}{\rm Tr \ Res }(dC)^{n-2}CdB\frac{dx}{x}=
{\rm Tr \ Res }C^{2n-3}dB\frac{dx}{x}={\rm Tr \ Res }2(dC)^{n-1}B\frac{dx}{x}.
\end{equation}
Thus
\begin{multline}{\rm Tr \ Res }(dC)^{n-1}B\frac{dx}{x}={\rm Tr \ Res
}(dC)^{n-2}
(CdB-dCB)\frac{dx}{x}=\\
 -{\rm Tr \ Res\ }d[C(dC)^{n-3}d(CB)\frac{dx}{x}].
\end{multline}
This form is exact, so we are done.
\end{proof}

\subsection{}
We want now to understand the image of $w_n(E, \nabla)$ under
the map $d_n$ defined in \ref{eqn5}. Of course \ref{prop2} says
that $d_n(w_n((E, \nabla)|(X-D))) = c_n(E)$.

\begin{defn}
(see \cite{EV}, Appendix B): Let $(E,
\nabla)$ be a flat connection with
logarithmic poles along $D$, with residue
$$
\Gamma = \oplus \Gamma_s \in \oplus_s H^0 (D_s , {\rm End}
E|_{D_s}). $$
One defines
\begin{gather}
N_i^{CH}(\Gamma) = (-1)^i \sum_{\alpha_1 + \cdots + \alpha_s=i}
\binom{i}{\alpha} {\rm Tr}(\Gamma_1^{\alpha_1} \circ \cdot ... \cdot
\circ \Gamma^{\alpha_s}_s) \cdot [D_1]^{\alpha_1} \cdots
[D_s]^{\alpha_s} \notag \\
\in CH^i(X) \otimes \C.
\end{gather}
One defines as usual the corresponding symmetric functions
$c_i^{CH}(\Gamma) \in CH^i(X) \otimes \C$ as polynomial with
$\Q$ coefficients in the Newton functions $N_i^{CH}(\Gamma)$. For example
\begin{gather}
c^{CH}_{2} (\Gamma) = \frac{1}{2} [(\sum_{s} {\rm Tr} \ (\Gamma_s) \cdot
D_s )^2 -  \notag \\
2 (\sum_{s} {\rm Tr} \ (\Gamma_s \cdot \Gamma_s ) \cdot
D^{2}_{s} + 2 \sum_{s <t} {\rm Tr} \ (\Gamma_s \cdot \Gamma_t ) D_s
\cdot D_t )]\\
\in CH^2(X) \otimes \C. \notag
\end{gather}
We denote by $c_2 (\Gamma)$ its image in $H^2 (X,
\Omega^{2}_{X,cl})$ and also by $c_2 (\Gamma)$ its image in $H^2
(X, \sH^{2}_{DR}$).
\end{defn}

Note that these invariants vanish when the connection has
nilpotent residues $\Gamma_s$. (This condition forces the local monodromies
around the components of $D$ to be unipotent (see \cite{D})).

\begin{thm}\label{Gamma}
Assume $k$ has characteristic zero and $X$ is proper.
Then
\begin{equation}{}
c_2 (E) - c_2 (\Gamma) = d_2 (w_2 (E, \nabla)) \in
H^2 (X, \sH^{2}).
\end{equation}
\end{thm}
\begin{proof}
In order to simplify the notations, we denote by $c_2
(\Gamma)$ the same expression in $CH^2 (X) \otimes \C$, $\oplus_s CH^1
(D_s) \otimes \C$, $\oplus_s F^1 H^{2}_{DR} (D_s)$ etc., where we always
distribute $2 D_s \cdot D_t$ for $s < t$ as one $D_s \cdot D_t$
on $D_s$ and one on $D_t$.

We denote by $\pi: Q \to X$ the flag bundle of $E$. As $\pi^*$
induces an isomorphism
\begin{equation}{}
\frac{H^2 (X, d \Omega^{1}_{X} )}{H^1 (X, \sH^{2}) } = \frac{H^3
(X, \sO_X \to \Omega^{1}_{X})}{N^1 H^{3} (X)} \> \sim >>
\frac{H^3 (Q, \sO_{Q} \to \Omega^{1}_{Q})}{N^1 H^{3} (Q)}
\end{equation}
and an injection $H^2 (X, \sH^{2}) \to H^2 (Q,
\sH^{2})$, it is enough to prove the compatibility on $Q$
via the exact sequence (\cite{BO})
\begin{equation}{}
0 \to \frac{H^3 (X, \sO_X \to \Omega^{1}_{X})}{N^1 H^{3}
(X)} \to H^2 (X, \Omega^{2}_{X, {\rm clsd}}) \to H^2 (X, \sH^{2}).
\end{equation}
Write $D'{}_s=\pi^*D_s$, and consider $(\sO(D'_s), \nabla_s)
\in \H^1 (Q, \sK_1 \to
\Omega^{1}_{Q} ({\rm log} D'_s)_{{\rm clsd}})$, where $\nabla_s$ is the
canonical connection with residue -1 along $D'_s$.

We define a product
\begin{gather}
(\sK^{m}_{i} \to (\pi^* \Omega^{i}_{X} ({\rm log} D))_{\tau d})
\times (\sK^{m}_{j} \to (\pi^* \Omega^{j}_{X} ({\rm log}
D))_{\tau d}) \notag
\\
\> \bullet >> (\sK^{m}_{i+j} \to (\pi^* \Omega^{i+j}_{X} ({\rm
log} D))_{\tau d})
\notag
\end{gather}
by
\begin{equation}x\cdot x' = \begin{cases}x\cup x' & \text{if deg }x'=0\\
           \tau d\log x\wedge x' & \text{if deg }x=0\text{ and deg}x'=1\\
0 & \text{ otherwise}
\end{cases}
\end{equation}
(Here $\tau d:\pi^*\Omega_X^i(\log D)\to
\pi^*\Omega_X^{i+1}(\log D)$ comes from the splitting
$\tau : \Omega^1_Q(\log D') \to \pi^*\Omega^1_X(\log D)$.
See proposition \ref{prop4.4.1}
as well as \cite{EI} and \cite{EII}.)
One verifies that
\begin{equation}{}
d (x \cdot x') = dx \cdot x' + (-1)^{{\rm deg} x} x \cdot dx',
\end{equation}
the only non trivial contribution left and right being for ${\rm
deg} x = {\rm deg} x' =0$.

This product defines elements ($W_1$ is the weight filtration)
\begin{equation}{}
\epsilon_{st} = ( \sO(D'_s), \nabla_s) \cdot (\sO (D'_t),
\nabla_t) \in \H^2 (Q, \sK_2 \to W_1 \Omega^{2}_{Q} ({\rm log}
(D'_s + D'_t))_{cl})
\end{equation}
which map to $D'_s \cdot D'_t$ in $CH^2 (Q)$. Moreover ${\rm Res\ }
\epsilon_{st}$ is the class of $D'_s \cdot D'_t$
sitting diagonally in
$$F^1 H^{2}_{DR}(D'_s) \oplus F^1 H^{2}_{DR} (D'_t)$$
if $s \neq t$; or in $F^1H^2_{DR} (D'_s)$ if $s =t$.

Next we want to define a cocycle $N_2(\pi^*(E, \nabla))$.

Let $h_{ij} (=h)$ be the upper triangular transition functions
of $E|_Q$ adapted to the tautological flag $E_i$, and write $B_i$ for the
local connection matrix in $ \Omega^{1}_{Q} ({\rm log} D'), D' = \pi^{-1} D$.
Then
$\tau B_i$ is upper triangular, and $\tau d B_i = d \tau B_i$
has zero's on the diagonal \cite{EI}, (0.7), (2.7). Let
$$w_i = \text{Tr}(B_idB_i).$$
Using ${\rm Tr} \ (d h h^{-1})^3
=0$, one computes that
 $w_i - w_j = - 3 {\rm Tr} \ d (h^{-1} dh B_j)$. But
\begin{gather}{}
{\rm Tr} \ h^{-1} dh B_j = {\rm Tr} \ h^{-1} B_i h B_j
\notag \\
{\rm Tr} \ h^{-1}_{ik} dh_{ij} dh_{jk} = {\rm Tr} \ h^{-1}_{ik} (B_i h_{ij} -
h_{ij} B_j) (B_j h_{jk} - h_{jk} B_j )
\\
= \delta {\rm Tr} \ (B_i h B_j h^{-1}) .\notag
\end{gather}
Here $\delta$ is the Cech coboundary. Writing $\mathcal C^i$ for Cech
$i$-cochains, we may define
\begin{gather}{}
3 N_2 (\pi^* (E, \nabla)) = ( 3 \sum^{r}_{a = 0} \xi^{a}_{ij}
\cup \xi^{a}_{jk} , - 3 {\rm Tr} \ (h^{-1} dh B_j), w_i) \in
\notag \\
(\sC^2 (Q, \sK_2) \times \sC^1 (Q, \Omega^{2}_{Q} ({\rm log} D')) \times
\sC^0 (Q, \Omega^{3}_{Q} ({\rm log} D') ))_{d +\delta}
\end{gather}
where $(\xi^{1}_{ij}, \ldots , \xi^{r}_{ij})$ is the
diagonal part of $h_{ij}$. This defines $3 N_2 (\pi^* (E,
\nabla))$ as a class in $\H^2 (Q, \sK_2 \to \Omega^{2}_{Q} ({\rm
log} D') \to ...)$ which maps to
\begin{gather}{}
3 \tau N_2 (\pi^* (E, \nabla)) = (3 \sum^{r}_{a=1}
\xi^{a}_{ij} \cup \xi^{a}_{jk}, 3 \sum^{r}_{a=1}
\omega^{a}_{i} \wedge (\delta \omega^a)_{ij}, 0)
\notag \\
\in \H^2 (Q, \sK_2 \to \pi^*  \Omega^{2}_{X} ({\rm log}
D')_{\tau d} )
\end{gather}
where $(\omega^{1}_{i}, \ldots ,
\omega^{r}_{i})$ is the diagonal part of $\tau B_i $.

As the image of $\tau N_2 (\pi^* (E, \nabla))$ in $H^2 (Q,
\sK_2)$ is just the second Newton class of $E$, the argument of
\cite{EII}, (1.7) shows that
\begin{gather}{}
N_2 (E, \nabla) : = \tau N_2 (\pi^* (E, \nabla)) \in \\
\H^2 (X,
\sK_2 \to \Omega^{2}_{X} ({\rm log} D) \to ...) \subset \H^2 (Q,
\sK_2 \to \pi^* \Omega^{2}_{X} ({\rm log} D) \to ...). \notag
\end{gather}

We
observe that $w(B) = {\rm Tr} \ B d B \in
W_2 \Omega^{3}_{Q} ({\rm log} D')$ (weight filtration) so the cocycle
\begin{gather}{}
2 x = - N_2 (\pi^* (E, \nabla)) + c_1 (\pi^* (E, \nabla))^2 =\\
\notag
(-{\rm Tr} \ (h^{-1} dh)^2 + {\rm Tr} \ h^{-1} dh \cdot {\rm Tr}
\ h^{-1} dh,\\
\notag
{\rm Tr} \ (h^{-1} dh B) - {\rm Tr} \ h^{-1} dh \cdot {\rm Tr} \
B,  - \frac{w(B)}{3}) \notag
\end{gather}
defines a class in
$$
\H^2 (Q, \Omega^{2}_{cl} \to W_1 \Omega^{2}_{Q} ({\rm log} D')
\to W_2 \Omega^{3}_{Q} ({\rm log} D')_{cl} ).
$$
One has an exact sequence
\begin{gather}{}
0 \to \H^2 (Q, \Omega^{2}_{cl} \to W_1 \Omega^{2}_{Q} ({\rm
log} D') \to W_1 \Omega^{3}_{Q} ({\rm log} D')_{cl})
\notag \\
\to \H^2 (Q, \Omega^{2}_{cl} \to W_1 \Omega^{2}_{Q} ({\rm log}
D') \to W_2 \Omega^{3}_{Q} ({\rm log} D')_{cl} )
\\
\> {\rm residue} >> \oplus_{s < t} H^0 (D'_{st},
\Omega^{1}_{D'_{st,cl}}). \notag
\end{gather}
As $D'_{st}$ is proper smooth, one has
$$
H^0 (D'_{st}, \Omega^{1}_{D'_{st,cl}}) \subset H^0 (D'_{st},
\sH^{1}) = H^{1} (D'_{st}).
$$
The residue of $2x$ along $D'_{st}$ is just the residue of $-
\frac{1}{3} w (B)$ along $D'_{st}$ via
\begin{gather}{}
H^0 (Q, \sH^{3} ({\rm log} D')) = H^0 (Q - D' ,
\sH^{3})
\notag \\
\to \oplus_s H^0 (D'_s -\cup_{t \neq s} D'_t, \sH^{2})
\\
\to \oplus_{s < t} H^0 (D'_{st}, \sH^{1}), \notag
\end{gather}
which vanishes.
Therefore
\begin{equation}{}
2 x \in \H^2 (Q, \Omega^{2}_{cl} \to W_1 \Omega^{2}_{Q} ({\rm
log} D') \to W_1 \Omega^{3}_{Q} ({\rm log} D')_{cl} ).
\end{equation}
Its residue in $\oplus_s H^1 (D'_s , \Omega^{1}_{D'_s})$ is
$({\rm Tr} \
(h^{-1} dh \cdot \Gamma) - {\rm Tr} \ h^{-1} dh \cdot {\rm Tr} \ \Gamma)$. By
\cite{EV}, Appendix B, one has
$- h^{-1} dh = \sigma (D') \cdot \Gamma$ in $H^1 (Q,
\Omega^{1}_Q{\otimes} {\rm End} E)$ where $\sigma (D')$ is the
extension
$$
0 \to \Omega^{1}_{Q} \to \Omega^{1}_{Q} ({\rm log} D') \to
\oplus_s \sO_{D'_s} \to 0.
$$
One has
\begin{enumerate}
\item $ - D'_s \cdot D'_s$ is the push down extension of
$\sigma (D'_s)$ by $\Omega^{1}_{Q} \to \Omega^{1}_{D'_s}$ in
$H^1 (Q, \Omega^{1}_{D'_s})$
\item  $ - D'_s \cdot D'_t$ is the extension
$$
0 \to \Omega^{1}_{D'_t} \to \Omega^{1}_{D'_t} ({\rm log} (D'_s
\cap D'_t)) \to \sO_{D'_s \cap D'_t} \to 0
$$
in $H^1 (Q, \Omega^{1}_{D'_t})$.
\end{enumerate}
It follows that residue $x = c_2 (\Gamma)$ in $\oplus_s H^1 (D'_s ,
\Omega^{1}_{D'_s})$. \\

For appropriate $\lambda_{st} \in k$ (the
coefficients of $c_2 (\Gamma)$), $c_2 (\Gamma) =$ residue $\sum
\lambda_{st} \epsilon_{st}$ in $\oplus_s H^1(D'_s, \Omega^1_{D'_s})$.
So one has
\begin{equation}
{\rm residue} (x - \sum \lambda_{st} \epsilon_{st}) \in
\oplus_s F^2 H^{2} (D'_s)
\end{equation}

Again, since residue $(x- \sum
\lambda_{st} \epsilon_{st}) = $ residue $x =0$ in $\oplus_s H^0
(D'_s , \sH^{2}) \subset \oplus_s H^{2} (k (D'_s))$,
one has that in $\oplus_s F^1 H^{2} (D'_s)$
\begin{equation}
{\rm residue} (x -
\sum \lambda_{st} \epsilon_{st} ) \in \oplus_s F^2 H^{2}
(D'_s) \cap H^1 (D'_s , \sH^{1}) = 0 .
\end{equation}

This shows that
residue $(x - \sum \lambda_{st} \epsilon_{st} ) =0$ in $\oplus_s
F^1 H^{2} (D'_s)$, that is
\begin{gather}{}
w_2 (E, \nabla) = (x - \sum \lambda_{st} \epsilon_{st}) \in
\notag \\
\frac{\H^2 (Q, \Omega^{2}_{cl} \to \Omega^2 \to \Omega^{3}_{cl})
= H^0 (Q, \sH^{3})}{{\rm Im} \oplus_s H^{1} (D'_s)}
\end{gather}
and maps to
\begin{equation}
c_2 (E) - c_2 (\Gamma) \ \mbox{in} \ H^2 (Q, \sH^{2}).
\end{equation}
\end{proof}

\subsection{question}\label{quest} We know (see \cite{EV}, Appendix B) that
on $X$ proper,
the image of $c_n(\Gamma)$ in the de Rham cohomology
$H^{2n}_{DR}(X)$ is the Chern class $c_n^{DR}(E)$.
This inclines to ask whether
$$c_n(E) - c_n(\Gamma) = d_n (w_n(E,
\nabla)) \in {\rm Griff}^n(X) .$$

\subsection{}
\begin{thm}\label{VHS} Let $(E, \nabla)$ be a flat connection with
logarithmic poles along a normal crossing divisor $D$
on a smooth proper variety $X$ over $\C$. When  $(E, \nabla)|(X-D)$
is a complex variation of Hodge structure, then $w_2(E, \nabla) = 0$
if and only if $c_2(E) - c_2(\Gamma) = 0 \in H^2(X, \sH^2)$. When
furthermore $(E, \nabla)|(X-D)$ is a Gau{\ss}-Manin system, then
$w_2(E, \nabla) \in H^0(X, \sH^3(\Q(2)))$, and if it has nilpotent
residues along the components of $D$, then
$w_2(E, \nabla) = 0$ if and only if  $c_2(E) = 0 \in H^2(X, \sH^2).$
\end{thm}
\begin{proof}
As in proposition \ref{prop3}, one has $w_n(E, \nabla) \in F^0$. In fact,
the proof does not use that $\nabla$ is everywhere regular, but only
that $w_n(E, \nabla)$ comes from a class in
$\H^n(Y, \sK^m_n \to \Omega^n_Y({\rm log}(Y-U) \to ...)$ on
some $(U, Y, \pi) \in \sS$. Further, if $(E, \nabla)|(X-D)$
is a $\C$ variation of Hodge structure,  then $w_n(E, \nabla) \in H^0(X,
\sH^{2n-1}(\R(n))$ as  $$w_n(E, \nabla)|(X-D) \in H^0(X-D,
\sH^{2n-1}(\R(n)))$$ (see proof of proposition 5.6.1).
When $(E, \nabla)|(X-D)$ is a Gau{\ss}-Manin system, then again one
argues exactly as in the proof of proposition 5.6.1 using \cite{CEs}
to show $w_2(E, \nabla) \in H^0(X, \sH^3(\Q(2)))$. Finally $c_2(\Gamma) =0$
when the residues of the connection are nilpotent.
\end{proof}

\section{Examples}\label{sec:ex}

\subsection{}\label{curves}
Let $X$ be a good compactification of the moduli space of
curves of genus $g$ with some level, such that a universal
family $\varphi : \sC \to X$ exists. Let $(E, \nabla)$ be the
Gau{\ss}-Manin system $R^1 \varphi_* \Omega^{\bullet}_{\sC/X} ({\rm
log} \infty)$. Then Mumford \cite{Mu}, (5.3) shows that
$c^{CH}_{i} (E) \otimes \Q = 0$ in $CH^i (X) \otimes \Q$ for $i
\geq 1$, so a
fortiori $c_i (E) \otimes \Q =0$ in the Griffiths group. As
$\nabla$ has nilpotent residues (the local monodromies at
infinity of the Gau{\ss}-Manin system are unipotent and $(E,
\nabla)$ is Deligne's extension \cite{D}), one applies
theorem \ref{VHS}.

In particular, for any semi-stable family of curves $\varphi : Y
\to X$ over a field $k$ of characteristic $0$,
$$
w_n (R^1 \varphi_* \Omega^{\bullet}_{Y/X} ({\rm
log} \infty)) = 0,
$$
for $n=2$ and for $n \ge 2$ if $\varphi$ is smooth (or if the question
(\ref{quest}) has a positive response).
\subsection{}
Let $X$ be a level cover of the moduli space of abelian
varieties such that a universal family $\varphi : \sA \to X$
exists. The Riemann-Roch-Grothendieck theorem applied to a
principal polarization $L$ on $\sA$ together with Mumford's
theorem that
$$
\varphi_* L^n = M \otimes \ \mbox{trivial}
$$
for some rank 1
bundle $M$, imply that $c^{CH}_{i}(E) \otimes \Q =0$ for $E =
R^1 \varphi_* \Omega^{\bullet}_{\sA/X}|{X_0}$, where
$X_0$ is the smooth locus of $\varphi$. This result was
communicated to us by G. van der Geer \cite{Ge}. In
particular, for any smooth family $\varphi : Y \to X$ of abelian
varieties with $X$ smooth proper over a field of characteristic $0$,
$w_n (R^1 \varphi_*
\Omega^{\bullet}_{Y/X} ) = 0$ for all $n \geq 2$.
\subsection{}
Let $\varphi: Y \to X$ be a smooth proper family of surfaces
over $X$ smooth. The Riemann-Roch-Grothendieck theorem, as
applied by Mumford in \cite{Mu}, implies that
the Chern character verifies
$$ {\rm ch} (\sum_{i=0}^{i=4}(-1)^i R^i\varphi_*
\Omega^{\bullet}_{Y/X}) \in CH^0(X) \otimes \Q \subset
CH^{\bullet}(X) \otimes \Q.$$ As $R^1\varphi_*
\Omega^{\bullet}_{Y/X}$ is dual to
$R^3\varphi_*\Omega^{\bullet}_{Y/X}$, the two previous examples
show that $c_i(R^2\varphi_*\Omega^{\bullet}_{Y/X}) = 0 $
in $CH^i(X) \otimes \Q$ for $i \geq 1$. This implies
$w_n (R^2 \varphi_*
\Omega^{\bullet}_{Y/X} ) = 0$ for all $n \geq 2$ when $X$ is proper.
\subsection{}
As shown in \cite{ECon}, $w_n (E, \nabla) =0$ in
characteristic zero when $(E, \nabla)$ trivializes on a finite
(not necessarly smooth) covering of $X$.
\subsection{}
Let $\varphi : Y \to X$ be a smooth proper family defined
over a perfect field $k$ of sufficiently large characteristic. Let $(E,
\nabla)$ be the Gau{\ss}-Manin system $R^a \varphi_*
\Omega^{\bullet}_{Y/X}$. Consider $w_n (E, \nabla) \in H^0 (X,
\sH^{2n-1})$, which is the restriction of the corresponding
class in characteristic zero when $\varphi$ comes from a smooth
proper family in characteristic zero. Assume this.
As $H^0(X, \sH^{2n-1}) \subset H^0(k(X), \sH^{2n-1})$, we
may assume that $R^a\varphi_*\Omega^{\bullet}_{Y/X}$ is locally
free and compatible with base change.

Via the diagram
\begin{equation}
\begin{CDS}
Y  \> F_{{\rm rel}} >> Y^{(p)} \>>> Y \\
\novarr \SE \varphi EE \V V \varphi^{(p)} V \novarr \V V \varphi V \\
\noharr X \> F >> X
\end{CDS}
\end{equation}
where $F$ is the absolute Frobenius of $X$, $\varphi^{(p)} =
\varphi \times_X F$, $F_{{\rm rel}}$ is the relative Frobenius,
one knows
by \cite{Ka}, (7.4) that the Gau{\ss}-Manin system $R^a \varphi_*
\Omega^{\bullet}_{Y/X} $ has a Gau{\ss}-Manin stable filtration
$R^a \varphi^{(p)}_{*} F_{{\rm rel}} (\tau_{\leq a}
\Omega^{\bullet}_{Y/X})$, such that the restriction of $\nabla$
to the graded pieces $F^* R^{a-i} \varphi_* \Omega^{i}_{Y/X}$ is
the trivial connection.

In particular, the graded pieces are locally generated by flat
sections and $A_i =0$ . So by additivity of the classes
$c_i (E, \nabla)$, $w _n(E, \nabla) = \theta_n(E, \nabla)= 0$.

In particular, the classes $w_n($ Gau{\ss}-Manin$)$ provide
examples of classes $w \in H^0 (X,
\sH^{2n-1})$ whose restriction modulo $p$ vanish for all but
finitely many $p$. This should imply that $w=0$ according
to \cite{Og}.

\bibliographystyle{plain}
\renewcommand\refname{References}

\end{document}